\documentclass[10pt]{article} % For LaTeX2e
\usepackage[utf8]{inputenc}
\usepackage[preprint]{tmlr}

\usepackage{amsmath,amsfonts,bm}

\def\1{\bm{1}}

\def\rmF{{\mathbf{F}}}
\def\rmG{{\mathbf{G}}}

\def\rmM{{\mathbf{M}}}

\def\rmO{{\mathbf{O}}}

\def\rmR{{\mathbf{R}}}

\def\rmV{{\mathbf{V}}}

\def\rmX{{\mathbf{X}}}
\def\rmY{{\mathbf{Y}}}
\def\rmZ{{\mathbf{Z}}}

\DeclareMathAlphabet{\mathsfit}{\encodingdefault}{\sfdefault}{m}{sl}
\SetMathAlphabet{\mathsfit}{bold}{\encodingdefault}{\sfdefault}{bx}{n}

\newcommand{\E}{\mathbb{E}}

\newcommand{\R}{\mathbb{R}}

\usepackage{url}

\usepackage{amssymb,enumerate}
\usepackage{enumitem} 
\usepackage{amsmath,amsfonts}
\usepackage{amsthm}
\usepackage{booktabs}
\usepackage{latexsym}
\usepackage{mathrsfs}
\usepackage{color}
\usepackage{microtype}
\usepackage{subfig}
\usepackage{makecell}%

\usepackage{algorithm}
\usepackage{algorithmicx}
\usepackage{algpseudocode}
\usepackage{graphicx}
\usepackage{multirow}
\allowdisplaybreaks[4]
\usepackage{xcolor}

\theoremstyle{plain}
\newtheorem{theorem}{Theorem}

\newtheorem{lemma}{Lemma}

\newtheorem{assumption}{Assumption}

\newtheorem{remark}{Remark}

\title{DeMuon: A Decentralized Muon for Matrix Optimization \\ 
over Graphs}

\author{\name Chuan He \email chuan.he@liu.se \\
	\addr Department of Mathematics\\
	Link\"oping University
	\AND
	\name Shuyi Ren \email shuyi.ren@liu.se\\
	\addr Department of Electrical Engineering\\
	Link\"oping University
	\AND
	\name Jingwei Mao \email jingwei.mao@liu.se\\
	\addr Department of Computer and Information Science\\
	Link\"oping University
	\AND
	\name Erik G. Larsson \email erik.g.larsson@liu.se\\
	\addr Department of Electrical Engineering\\
	Link\"oping University}

\def\month{MM}  % Insert correct month for camera-ready version
\def\year{YYYY} % Insert correct year for camera-ready version
\def\openreview{\url{https://openreview.net/forum?id=XXXX}} % Insert correct link to OpenReview for camera-ready version

\begin{document}

\maketitle

\begin{abstract}
%	In this paper, we propose {\it DeMuon}, a new method for decentralized matrix optimization over a communication graph. {\it DeMuon} incorporates matrix orthogonalization via Newton-Schulz iterations---a technique inherited from its centralized predecessor, {\it Muon}---and employs gradient tracking to mitigate heterogeneity among local functions. Under suitable assumptions, we establish its convergence guarantees in terms of stationarity and consensus errors. To the best of our knowledge, {\it DeMuon} is the first direct extension of {\it Muon} to decentralized optimization over graphs with provable complexity guarantees. Finally, we conduct numerical experiments on decentralized transformer pretraining over graphs with varying degrees of connectivity.  Our numerical results demonstrate a clear margin of improvement of {\it DeMuon}
%	over other popular decentralized algorithms.

This paper considers fully decentralized matrix optimization over communication graphs. We propose {\it DeMuon}, a decentralized matrix optimization method that incorporates matrix orthogonalization, a technique inherited from its centralized predecessor, {\it Muon}, and employs gradient tracking to mitigate heterogeneity among local objective functions. Under spectral- and nuclear-norm-based gradient Lipschitz assumptions, we establish global convergence guarantees for {\it DeMuon} in terms of stationarity and consensus errors. Furthermore, we propose an accelerated variant, named {\it DeMuon-A}, which integrates the multi-extrapolation techniques. Under an additional higher-order Lipschitz assumption induced by the spectral and nuclear norms, we establish improved convergence guarantees for {\it DeMuon-A}. Finally, numerical experiments demonstrate the practical advantages of the proposed methods. 

%    {\color{blue}I will rewrite the abstract.}
\end{abstract}


\section{Introduction}\label{sec:intro}

Recently, a matrix optimization method, called {\it Muon} \citep{jordan2024muon}, has attracted significant attention from both the deep learning and optimization communities. {\it Muon} departs from traditional methods, such as {\it Adam} \citep{kingma2015adam}, which are typically based on vectorizing matrix variables, and demonstrates advantages through strong empirical results on the training of language models. At a high level, when applied to the matrix optimization problem $\min_{X\in\R^{m\times n}} f(X)$, {\it Muon} generates a sequence $\{X^{(k)}\}$ according to
\begin{align*}
M^{(k)} = (1-\theta) M^{(k-1)} + \theta G(X^{(k)};\xi^{(k)}),\quad X^{(k+1)} = \underset{\|X-X^{(k)}\|\le\eta}{\mathrm{arg\,min}}\ \langle M^{(k)},X-X^{(k)}\rangle\qquad\forall k=0,1,\ldots,
\end{align*}
where $G:\R^{m\times n}\times\Xi\to\R^{m\times n}$ is the stochastic gradient estimator %(see Assumption \ref{asp:basic}(c) for precise conditions), 
and $\|\cdot\|$ denotes the spectral norm of a matrix. Looking at the update, {\it Muon} leverages the momentum update for accumulating stochastic gradient estimates (common in deep learning optimization), together with a norm-constrained linear subproblem that serves as a local search oracle (which is quite popular, as it enforces normalization). {\it Muon}'s success is attributed to its use of the spectral norm in constructing linear subproblems, which appears particularly well suited to a variety of neural network architectures.

% Muon is a simple algorithm, combine momentum and linear minimization subrpoblem over spectral norm ball. (one sentence, explain the update) The core magic/effective is the use of spectal norm replacing previous Frobenious norm when dealing with matrix optimizaiton problems. (core advantage)

Beyond its widely recognized empirical advantages in centralized deep learning optimization \citep{jakovetic2014fast,liu2025muon}, {\it Muon} has also demonstrated strong performance gains on many classical centralized matrix optimization problems, including matrix regression and matrix completion (e.g., see \cite{lau2025polargrad,he2025low}). Moving to the decentralized world, matrix optimization problems are likewise ubiquitous, arising in applications such as distributed deep learning \citep{nazari2022dadam} and classical problems such as distributed PCA \citep{liang2013distributed} and distributed matrix completion and factorization \citep{mackey2015distributed}. Motivated by the prevalence of such applications, in this paper, we take a step toward bringing {\it Muon}'s techniques to the decentralized world, where we formulate the following multi-agent optimization problem:
\begin{align}\label{pb:finite-sum-matrix}
\min_{X\in\R^{m\times n}} f(X) := \frac{1}{N}\sum_{i=1}^N f_i(X),
\end{align}
where $f_i:\R^{m\times n}\to\R$ denotes the local objective associated with node $i$ for $i\in[N]$. The $N$ agents are connected through a graph $\mathcal{G}=(\mathcal{V},\mathcal{E})$, where $\mathcal{V}=[N]$ denotes the set of all node indices and $\mathcal{E}$ consists of all directed pairs $(i,j)\in[N]\times [N]$ such that node $i$ can send information to node $j$.

% .... Muon is successful in matrix optimization, not just for deep learning [?] but also has been observed having advantages for matrix machine learning problem like matrix factorization cite (Polargrad. Low rank Orth.). % These problems are also wide spread in decentralized setting. Decnetralized deep learnig, Decentralized matrix machine learnign, 
% Motivated by these, we consider a general form of decentralized matrix opt ...

The contributions of this paper are threefold.
\begin{itemize}[leftmargin = 1.5em]
    \item We propose and analyze a fully decentralized variant of {\it Muon}, named {\it DeMuon}, for solving problem \eqref{pb:finite-sum-matrix}, which integrates spectral-norm-based updates into a gradient-tracking-based decentralized algorithmic framework. Under the spectral- and nuclear-norm-based gradient Lipschitz assumption, we establish the global convergence guarantees of {\it DeMuon} in terms of stationarity and concensus errors.
    %\textcolor{red}{add acceleration part}
    \item We propose and analyze an accelerated variant of {\it DeMuon}, named {\it DeMuon-A}, extending the multi-extrapolation techniques developed for vector-variate optimization problems \citep{he2025heavytail}. Under the higher-order Lipschitz assumption induced by spectral- and nuclear-norms, we establish the improved global convergence guarantees of {\it DeMuon-A} in terms of stationarity and concensus errors.
    
    %From a theoretical perspective, we establish a concentration inequality based approach to derive the convergence complexity of the proposed method. Notably, our bounds are dimension-free, in contrast to most existing works, which typically rely on standard matrix norm inequalities and yield dimension dependent guarantees under the Frobenius norm. 
    \item We compare the empirical performance of our {\it DeMuon} and {\it DeMuon-A} against competing methods on distributed deep learning and distributed matrix regression tasks. Numerical results demonstrate the clear practical advantages of {\it Muon}-type methods in decentralized settings.
    
    %From an empirical perspective, we demonstrate the effectiveness of our approach on both decentralized deep learning and distributed matrix optimization tasks.
\end{itemize}

%Theoretically, we derive contraction inequality based approach to derive complexity, our bound is dimension free comparing with most of existing works which are mainly proved based on the 

We now review recent algorithmic developments in decentralized settings that are most closely related to our study. A more complete review of related work is deferred to Appendix~\ref{Appendix:related_works}. 

\cite{takezawa2025fedmuon} proposed {\it FedMuon}, which generalizes {\it Muon} to the federated learning setting, where a central server coordinates communication and aggregates local updates from participating clients. However, this differs from our study, as we consider a fully decentralized setting without a central server.

In the fully decentralized setting without a central server, \cite{yu2026decentralized} proposed a decentralized gradient-tracking-based {\it SGD} method with normalization and momentum. If we apply their method to problem \eqref{pb:finite-sum-matrix}, the normalization can be interpreted as constructing a linear subproblem with a Frobenius norm constraint. The relationship between their method and {\it DeMuon} in the fully decentralized setting is analogous to the relationship between normalized {\it SGD} with momentum and {\it Muon} in the centralized setting.

%and show that directly using Muon in FedAvg fails to converge due to the bias of the LMO operator. To address this, they propose FEDMUON and show that it converges for any number of Newton--Schulz iterations, with more accurate LMO solutions yielding faster convergence. 

%\cite{yu2026decentralized} show that DSGD with gradient tracking, momentum, and normalization achieves optimal convergence rates under heavy-tailed noise. Inspired by their work, we released the first version of DeMuon last October and developed a new analytical framework establishing global convergence guarantees for DeMuon, which incorporates spectral-norm-based updates, gradient tracking and momentum.

After releasing our first version in October 2025, \cite{zhang2026suda} recently developed {\it SUDA-Muon}, which generalizes {\it DeMuon} by incorporating more advanced techniques for consensus optimization. Our work differs from theirs in that we adopt matrix martingale moment inequalities (see Lemma \ref{lem:concentration}) to establish convergence for {\it DeMuon} without explicit dependence on the problem dimensionality $(m,n)$. Moreover, we incorporate extrapolation techniques into {\it DeMuon} to achieve acceleration under higher-order smoothness assumptions. %Combining {\it Muon} with other advanced fully decentralized frameworks appears to be a promising direction, as suggested by the numerical results in \cite{zhang2026suda}. This line of research is orthogonal to the development of our paper, and our techniques can be readily integrated with theirs.

%Their empirical results demonstrate the advantages of fully decentralized Muon in decentralized consensus optimization.

%which allows ED/D$^2$, EXTRA, and ATC-style gradient tracking to serve as modular backbone choices. Their work generalizes our proposed method by incorporating several recent decentralized consensus optimization algorithms and provides corresponding theoretical convergence guarantees. The results clearly demonstrate the advantages of fully decentralized Muon in decentralized consensus optimization. Their focus is broader and largely orthogonal to ours. In contrast, we concentrate on the fundamental problem of establishing tight theoretical guarantees for our specific decentralized Muon variant through matrix martingale moment inequalities. Moreover, under higher-order Lipschitz assumptions, we further incorporate our recent acceleration technique into decentralized Muon, which, to the best of our knowledge, is the first work to do so.

%{\color{blue}discuss FedMuon, Distributed Normalized SGD, SUDA-Muon here.}

The remainder of this paper is organized as follows. In Section \ref{sec:not-ass}, we introduce the notation and assumptions used in the paper. In Section \ref{sec:d-muon}, we propose and analyze {\it DeMuon}. In Section \ref{sec:DeMuon+}, we propose and analyze {\it DeMuon-A}. Section \ref{sec:sr} presents the numerical experiments. In Appendix \ref{Appendix:related_works}, we include a review on relevant literature on {\it Muon} and decentralized optimization. Appendices \ref{sec:pf}-\ref{apx:pf-DeMuonA} contain the proofs of our main results.

\section{Notation and assumptions}
\label{sec:not-ass}
Throughout this paper, we use $\mathbb{R}^{m\times n}$ to denote the Euclidean space of $m\times n$ real matrices. We use $\|\cdot\|$ and $\|\cdot\|_*$ to denote the spectral norm and the nuclear norm of a matrix, respectively. We use $\langle\cdot,\cdot\rangle$ to denote the trace inner product for matrices. For any positive integer $p$ and a $p$th-order continuously differentiable function $\varphi:\R^{m\times n}\to\R$, we denote by $\mathcal{D}^p\varphi(X)[H_1,\ldots,H_p]$ the $p$th-order directional derivative of $\varphi$ at $X$ along $H_i\in\R^{m\times n}$, $i\in[p]$, and by $\mathcal{D}^p\varphi(X)[\cdot]$ the associated symmetric $p$-linear form. For any symmetric $p$-linear form $\mathcal{T}[\cdot]$, we define its operator norm induced by the spectral norm as
\begin{align}\label{def:operator-norm}
    \|\mathcal{T}\|_* = \max_{H_1,\ldots,H_p} \{\mathcal{T}[H_1,\ldots,H_p] : \|H_i\|\le1,\ i\in[p]\}.
\end{align}
Here, we slightly abuse notation, as for $p = 1$, this operator norm coincides with the nuclear norm. For any $X\in\R^{m\times n}$ and $H_i\in\R^{m\times n}$ with $i\in[p-1]$, we define $\nabla^p\varphi(X)(H_1,\ldots,H_{p-1})\in\R^{m\times n}$ by
\begin{align*}
    \langle \nabla^p\varphi(X)(H_1,\ldots,H_{p-1}), H_p\rangle := \mathcal{D}^p\varphi(X)[H_1,\ldots,H_p]\qquad\forall H_p\in\R^{m\times n}.
\end{align*}
For any $X,H\in\R^{m\times n}$, we let $\mathcal{D}^p\varphi(X)[H]^p:=\mathcal{D}^p\varphi(X)[H,\ldots,H]$ and $\nabla^p\varphi(X)(H)^{p-1}:=\nabla^p\varphi(X)(H,\ldots,H)$.

We define the matrix sign of any nonzero $M \in \mathbb{R}^{m\times n}$ as $\mathrm{msgn}(M) = UV^T$, where $U \in \mathbb{R}^{m\times r}$ and $V \in \mathbb{R}^{n\times r}$ are orthonormal matrices obtained from the reduced SVD of $M$. For any two real symmetric
matrices $A_1$ and $A_2$, we write $A_1\preceq A_2$ (respectively, $A_1\prec A_2$) if $A_2-A_1$ is positive semidefinite (respectively, definite). For a collection of $m \times n$ matrices $\{X_i\}_{i \in [N]}$, we define the stacked notation:
\begin{align*}
\mathbf{X}:=\left[\begin{matrix}
		X_1\\
		\vdots\\
		X_N
	\end{matrix}\right],\quad \mathbf{F}(\mathbf{X}):=\left[\begin{matrix}
		f_1(X_1)\\
		\vdots\\
		f_N(X_N)
	\end{matrix}\right],\quad \nabla \rmF(\rmX):=\left[\begin{matrix}
		\nabla f_1(X_1)\\
		\vdots\\
		\nabla f_N(X_N)
	\end{matrix}\right],
\end{align*}
and define the averaged notation as
\begin{align*}
\Bar{X}:=\frac{1}{N}\sum_{i=1}^N X_i,\quad \bar{g}(\mathbf{X}):=\frac{1}{N}\sum_{i=1}^N\nabla f_i(X_i).   
\end{align*} 
We follow the convention that bold symbols denote variables in the stacked matrix space $\R^{(Nm)\times n}$, whereas italic mathematical symbols denote variable in $\R^{m\times n}$. Let $A\otimes B$ denote the Kronecker product. We let $I_d$ denote the $d\times d$ identity matrix for any integer $d\ge1$, $\mathbf{1}$ denote the $N$-dimensional all-ones vector, and $\mathbf{0}$ denote a all-zero matrix. In addition, we use $\mathcal{O}(\cdot)$ to denote the standard big-O notation.

We now make the following assumptions throughout this paper.

\begin{assumption}\label{asp:basic}
	\begin{enumerate}[leftmargin=2em]
		\item[{\rm (a)}] There exists $f_{\mathrm{low}}\in\R$ such that $f(X)\ge f_{\mathrm{low}}$ for all $X\in\R^{m\times n}$.
		\item[{\rm (b)}] There exists $L_*>0$ such that $\|\nabla f_i(X_i) - \nabla f_i(Y_i)\|_* \le L_{*}\| X_i -  Y_i\|$ for all $X_i,Y_i\in\R^{m\times n}$ and $i\in[N]$.
		\item[{\rm (c)}] There exists $V\in\R^{m\times n}$ such that stochastic gradient estimators $G_i:\R^{m\times n}\times\Xi\to\R^{m\times n}$ satisfy
		\begin{align*}
			\E[G_i(X_i;\xi_i)]=\nabla f_i(X_i),\quad \E[(G_i(X_i;\xi_i) - \nabla f_i(X_i))^T(G_i(X_i;\xi_i) - \nabla f_i(X_i))] \preceq V^TV\qquad    
		\end{align*}
		for all $X_i\in\R^{m\times n}$ and each $i\in[N]$.
		\item[{\rm (d)}] The mixing matrix $W\in\R^{N\times N}$ associated with the graph $\mathcal{G}$ has the properties: 
		\begin{enumerate}[leftmargin=1.6em]
			\item[{\rm (i)}] Primitivity: $W\ge0$, and $W^j>0$ for some positive integer $j$;
			\vspace{1mm}
			\item[{\rm (ii)}] Doubly stochasticity: $\mathbf{1}^TW=\mathbf{1}^T$, and $W\mathbf{1} = \mathbf{1}$.
            \vspace{1mm}
		\end{enumerate}
	\end{enumerate}    
\end{assumption}

\begin{remark}
	\begin{enumerate}[leftmargin=2em]
		\item[{\rm (i)}] Assumption \ref{asp:basic}(b) is common in the study of Muon and its variants (e.g., see \cite{he2025low,shen2025convergence}). Since $f$ is the average of all $f_i$'s, Assumption \ref{asp:basic} implies that $\|\nabla f(X)-\nabla f(Y)\|_* \le L_*\|X-Y\|$ for all $X,Y\in\R^{m\times n}$. It implies the following descent inequality (e.g., see \cite{shen2025convergence}):
		\begin{align}\label{ineq:desc}
			f(Y)\le f(X) + \langle\nabla f(X), Y-X\rangle + \frac{L_*}{2} \|Y-X\|^2\qquad\forall X,Y\in\R^{m\times n}.
		\end{align}
		\item[{\rm (ii)}] Assumption \ref{asp:basic}(c) is a bounded noise condition on each gradient estimator $G_i$, which has been adopted for fine-grained analyses of matrix-variate algorithms (e.g., see \cite{an2025asgo,pan2025unbiased}). It implies the widely used bounded variance condition $\E[\|G_i(X_i;\xi_i) - \nabla f_i(X_i)\|_F^2]\le\|V\|_F^2$.  
		\item[{\rm (iii)}] Under Assumption~\ref{asp:basic}(d), we define the rate of the mixing matrix $W$ as
		\begin{align}\label{upbd:eig-W}
			\lambda := \Big\|W - \frac{1}{N}\mathbf{1} \mathbf{1}^T \Big\| \in(0,1),
		\end{align}
		which is widely used to characterize the consensus performance of $W$ (see, e.g., \cite{xiao2004fast}).
	\end{enumerate}
\end{remark}

We adopt the following notation, which will be used frequently later:
        \begin{align}\label{def:Llambda-objgap}
		\Delta_f := f(\bar{X}^{(0)}) - f_{\mathrm{low}},\quad L_\lambda  := \bigg(\frac{2\sqrt{N}\lambda}{1-\lambda} + 1\bigg)L_*.
	\end{align}

\section{DeMuon: A decentralized Muon}\label{sec:d-muon}

In this section, we propose a decentralized variant of {\it Muon}, called {\it DeMuon} for brevity, and establish upper bounds on its consensus and stationarity errors. We also provide convergence guarantees for these errors when the algorithm parameters are properly selected.

\begin{algorithm}[!htbp] 
	\caption{DeMuon: A decentralized Muon} 
	\label{alg:r-msgn-1} 
	\begin{algorithmic}[0] 
		\State \textbf{Input:} starting iterate $\rmX^{(0)}$ with $X^{(0)}_i=X^{(0)}_j$ for any $i,j\in[N]$, mixing matrix $W\in\R^{N\times N}$, step length $\eta>0$, momentum parameter $\theta\in(0,1)$.
		\State \textbf{Initialize:} $\rmM^{(-1)}=\rmV^{(-1)}=\mathbf{0}$. 
		\For{$k=0,1,2,\ldots$} 
		\State Update the local gradient estimators:
		\begin{align}\label{update-mk} 
			M^{(k)}_i = (1-\theta)M^{(k-1)}_i +\theta  G_i(X^{(k)}_i;\xi^{(k)}_i)\qquad\forall i\in[N].  
		\end{align} 
		\State Update the global gradient estimators: 
		\begin{align}\label{update-vk} 
			V^{(k)}_i = \sum_{j=1}^N w_{ij} (V^{(k-1)}_j + M^{(k)}_j - M^{(k-1)}_j)\qquad\forall i\in[N]. 
		\end{align}
		\State Update the local iterates:
		\begin{align}\label{update-xk}
			X^{(k+1)}_i = \sum_{j=1}^N w_{ij} (X^{(k)}_j - \eta\cdot\mathrm{msgn}(V^{(k)}_j))\qquad\forall i\in[N].
		\end{align} 
		\EndFor
	\end{algorithmic} 
\end{algorithm}

Our proposed {\it DeMuon} generates three sequences, $\{\rmM^{(k)}\}$, $\{\rmV^{(k)}\}$, and $\{\mathbf{X}^{(k)}\}$. Specifically, at each iteration $k\ge0$, each node $i$ updates the local gradient estimator $M^{(k)}_i$ via an exponentially weighted moving average of the stochastic gradient $G_i$ evaluated at $X^{(0)}_i,\ldots,X^{(k)}_i$. Next, it applies a tracking technique (see, e.g., \cite{di2016next}) to update $V_i^{(k)}$, which ensures that $V_i^{(k)}$'s achieve consensus and approximate the global gradient. Then, the local iterate $X^{(k+1)}_i$ is updated by aggregating the orthogonalization updates performed on $V^{(k)}_j$'s from neighboring nodes. Details of {\it DeMuon} are described in Algorithm \ref{alg:r-msgn-1}.

The following theorem provides an upper bound on the consensus error of $\{\mathbf{X}^{(k)}\}$ generated by {\it DeMuon}. Its proof is deferred to Section \ref{sec:pf-thm1}.

\begin{theorem}[{{\bf consensus error}}]\label{thm:cs-error}
	Suppose that Assumption \ref{asp:basic} holds. Let $\{\mathbf{X}^{(k)}\}$ be generated by Algorithm \ref{alg:r-msgn-1} with step size $\eta>0$, and let $\lambda$ be defined in \eqref{upbd:eig-W}. Then it holds that for all $k\ge0$, 
	\begin{align}\label{ineq:cs-error}
		\|\mathbf{X}^{(k)} - \mathbf{1}\otimes\bar{X}^{(k)}\| \le \frac{\sqrt{N}\lambda\eta}{1-\lambda}.
	\end{align}
\end{theorem}

\begin{remark}
	Notice that the consensus error established in Theorem \ref{thm:cs-error} is similar to that for decentralized normalized vector-variate algorithms; see \cite[Eq.17]{yu2026decentralized}. On closer inspection, the consensus error in Theorem \ref{thm:cs-error} is measured using the spectral norm. Indeed, when applied to our matrix-variate optimization problem \eqref{pb:finite-sum-matrix}, the consensus error in Eq.17 in \citet{yu2026decentralized} reduces to 
	\begin{align*}
		\|\mathbf{X}^{(k)} - \mathbf{1}\otimes\bar{X}^{(k)}\|_F \le \frac{\sqrt{N}\lambda\eta}{1-\lambda},
	\end{align*}
	which is measured by the Frobenious norm. This mismatch is because normalized direction, i.e., $V_i^{(k)}/\|V_i^{(k)}\|_F$, lies on the Frobenius norm-induced unit sphere, whereas the matrix-sign direction, i.e., $\mathrm{msgn}(V_j^{(k)})$, lies on the spectral norm-induced unit sphere.
\end{remark}

We next derive a bound for the stationarity error for $\{\rmX^{(k)}\}$ generated by {\it DeMuon} in the following theorem, whose proof is deferred to Section \ref{sec:pf-thm2}.

\begin{theorem}[{{\bf stationarity error}}]\label{thm:stat-error}
Suppose that Assumption \ref{asp:basic} holds. Let $L_*$ and $V$ be given in Assumption \ref{asp:basic}, and let $\lambda$ and $(\Delta_f,L_\lambda)$ be defined in \eqref{upbd:eig-W} and \eqref{def:Llambda-objgap}, respectively. Let $\{\rmX^{(k)}\}$ be generated by Algorithm \ref{alg:r-msgn-1} with input parameters $\eta$ and $\theta$. Then,
	\begin{align}
		\frac{1}{K} \sum_{k=0}^{K-1} \E[\|\bar{g}(\rmX^{(k)})\|_*] \le \frac{\Delta_f}{K\eta} + \frac{4\sqrt{N\theta}\|V\|_*}{(1-\theta)(1-\lambda)} + \frac{(8N+1)L_\lambda \eta}{2\theta(1-\theta)(1-\lambda)} + \frac{2\theta\lambda\sqrt{N}\|V\|_*}{(1-\theta)\sqrt{1-\lambda}} + & \frac{4\|\nabla \rmF(\rmX^{(0)})\|_*}{K\theta(1-\theta)(1-\lambda)}\nonumber\\
		&\qquad\quad\forall K\ge1.\label{ineq:stat-error}
	\end{align}
\end{theorem}

\begin{remark}
	Observe that the right-hand side of \eqref{ineq:stat-error} exhibits a similar pattern in terms of dependence on the step size $\eta$ and the momentum parameter $\theta$ to that of the normalized SGD with momentum; see \cite[Theorem 1]{cutkosky2020momentum}. Taking $\eta=\mathcal{O}(K^{-3/4})$ and $\theta=\mathcal{O}(K^{-1/2})$ (optimal choices with respect to $K$), we can see that the average stationarity error is of the order of $\mathcal{O}(K^{-1/4})$, which matches the optimal convergence rate in the centralized setting. To derive the values of $\eta$ and $\theta$ for achieving a tight average stationarity error, we notice from \eqref{ineq:stat-error} that, with the choice of $\eta=\mathcal{O}(K^{-3/4})$ and $\theta=\mathcal{O}(K^{-1/2})$, the first three terms %$\frac{\Delta_f}{K\eta}$, $\frac{4\sqrt{N\theta}\|V\|_*}{(1-\theta)(1-\lambda)}$, and $\frac{(8N+1)L_\lambda \eta}{2\theta(1-\theta)(1-\lambda)}$ 
    dominate the dependence on $K$ among all terms on the right-hand side of \eqref{ineq:stat-error}. Therefore, we choose $\eta$ and $\theta$ to optimize the sum of these three terms, as done in the next theorem.
\end{remark}

The following theorem provides a convergence guarantee for {\it DeMuon} when $\eta$ and $\theta$ are selected properly. Its proof is relegated to Section \ref{sec:pf-thm3}.

\begin{theorem}[{{\bf convergence}}]\label{thm:conv}
Suppose that Assumption \ref{asp:basic} holds. Let $L_*$ and $V$ be given in Assumption \ref{asp:basic}, and let $\lambda$ and $(\Delta_f,L_\lambda)$ be defined in \eqref{upbd:eig-W} and \eqref{def:Llambda-objgap}, respectively. Let $K$ be the maximum iteration number of Algorithm \ref{alg:r-msgn-1} such that $K\ge\frac{4(1-\lambda)\Delta_fL_\lambda}{\|V\|_*^2}$, and define
	\begin{align}\label{def:hat-theta-hat-eta}
		\hat\theta = \frac{1}{\|V\|_*}\sqrt{\frac{(1-\lambda) \Delta_fL_\lambda }{K}},\quad \hat\eta = \sqrt{\frac{(1-\lambda)\Delta_f\hat\theta}{NL_\lambda K}}.
	\end{align}
	Let $\{\rmX^{(k)}\}$ be generated by Algorithm \ref{alg:r-msgn-1} with inputs $(\eta,\theta)=(\hat{\eta},\hat{\theta})$. Then, 
	\begin{align}
		%\frac{1}{K} \sum_{k=0}^{K-1} \E[\|\bar{g}(\rmX^{(k)})\|_*] \le \mathcal{O}\bigg(\frac{N^{\frac{5}{8}}\|V\|_*^{\frac{1}{2}}}{1-\lambda}\bigg(\frac{L_*}{K}\bigg)^{\frac{1}{4}}\bigg),\quad \|\rmX^{(k)} - \mathbf{1}\otimes\bar{X}^{(k)}\| \le \mathcal{O}\bigg(\frac{1}{N^{\frac{1}{8}}\|V\|_*^{\frac{1}{2}}L_*^{\frac{1}{4}}K^{\frac{3}{4}}}\bigg). \label{upbd:cs-ave}
        \frac{1}{K} \sum_{k=0}^{K-1} \E[\|\bar{g}(\rmX^{(k)})\|_*] \le \mathcal{O}(K^{-1/4}),\quad \max_{0\le k\le K-1}\{\|\rmX^{(k)} - \mathbf{1}\otimes\bar{X}^{(k)}\|\} \le \mathcal{O}(K^{-3/4}).\label{upbd:cs-ave}
	\end{align}
\end{theorem}

\begin{remark}
	The dependence on K of the stationarity error in \eqref{upbd:cs-ave} matches the results for centralized normalized SGD with momentum in \cite[Theorem 1]{cutkosky2020momentum}, while the consensus error matches that for decentralized normalized SGD in \cite[Lemma 4]{yu2026decentralized}, with the tail index equal to $2$.
\end{remark}

\section{DeMuon-A: An accelerated decentralized Muon}\label{sec:DeMuon+}

In this section, we propose an accelerated variant of {\it DeMuon}, denoted {\it DeMuon-A}. We establish upper bounds on its consensus error and stationarity error. We also provide convergence guarantees for these errors when the algorithm parameters are properly chosen.

{\it DeMuon-A} extends the multi-extrapolation technique proposed in \cite{he2025heavytail} for the centralized setting, which exploits higher-order smoothness of the nonconvex objective function to achieve acceleration. Our proposed {\it DeMuon-A} generates four sequences, $\{\rmZ^{(k,s)}\}$, $\{\rmM^{(k)}\}$, $\{\rmV^{(k)}\}$, and $\{\mathbf{X}^{(k)}\}$. Specifically, at each iteration $k \ge 0$, each node $i$ performs $q$ separate local extrapolation steps based on $X_i^{(k)}$ and $X_i^{(k-1)}$ to obtain an extrapolated point $Z_i^{(k,s)}$, for $s \in [q]$. Then, each node $i$ updates the local gradient estimator $M_i^{(k)}$ via an exponentially weighted moving average of the stochastic gradient $G_i$ evaluated at the extrapolated points $\{Z_i^{(k,s)}\}_{0\le t\le k,1\le s\le q}$. Next, it applies a tracking technique to update $V_i^{(k)}$, which ensures that $V_i^{(k)}$'s achieve consensus and approximate the global gradient. Finally, the local iterate $X^{(k+1)}_i$ is updated by aggregating the orthogonalization updates performed on $V^{(k)}_j$'s from neighboring nodes. Details of {\it DeMuon-A} are described in Algorithm \ref{alg:r-msgn-1+}.

\begin{algorithm}[t] 
	\caption{DeMuon-A: An accelerated decentralized Muon} 
	\label{alg:r-msgn-1+} 
	\begin{algorithmic}[0] 
		\State \textbf{Input:} starting iterate $\rmX^{(0)}$ with $X^{(0)}_i=X^{(0)}_j$ for any $i,j\in[N]$, mixing matrix $W\in\R^{N\times N}$, step length $\eta>0$, extrapolation count $q\ge1$, extrapolation parameters $\{\gamma_s\}_{s\in[q]}\subset(0,1)$, momentum parameters $\{\theta_s\}_{s\in[q]}$ with $\sum_{s=1}^q\theta_s\in(0,1)$. 
		\State \textbf{Initialize:} $\rmM^{(-1)}=\rmV^{(-1)}=\mathbf{0}$. 
		\For{$k=0,1,2,\ldots$} 
		\State Perform the local extrapolations:
		\begin{align}\label{update-zk+} 
			Z^{(k,s)}_i = X_i^{(k)} + \frac{1 - \gamma_s}{\gamma_s}(X_i^{(k)} - X_i^{(k-1)})\qquad \forall i\in[N], s\in[q].
		\end{align}
		\State Update the local gradient estimators:
		\begin{align}\label{update-mk+} 
			M^{(k)}_i = \Big(1-\sum_{s=1}^q\theta_s\Big)M^{(k-1)}_i + \sum_{s=1}^q\theta_s  G_i(Z^{(k,s)}_i;\xi^{(k)}_i)\qquad\forall i\in[N].  
		\end{align} 
		\State Update the global gradient estimators: 
		\begin{align}\label{update-vk+} 
			V^{(k)}_i = \sum_{j=1}^N w_{ij} (V^{(k-1)}_j + M^{(k)}_j - M^{(k-1)}_j)\qquad\forall i\in[N]. 
		\end{align}
		\State Update the local iterates:
		\begin{align}\label{update-xk+}
			X^{(k+1)}_i = \sum_{j=1}^N w_{ij} (X^{(k)}_j - \eta\cdot\mathrm{msgn}(V^{(k)}_j))\qquad\forall i\in[N].
		\end{align} 
		\EndFor
	\end{algorithmic} 
\end{algorithm}

To adapt to {\it Muon}-type algorithms, we impose a $(\|\cdot\|,\|\cdot\|_*)$-induced higher-order Lipschitz continuity condition, under which {\it DeMuon-A} achieves acceleration. This is analogous to the $(\|\cdot\|,\|\cdot\|_*)$-induced gradient Lipschitz condition in Assumption \ref{asp:basic}(b).

\begin{assumption}\label{asp:high-order}
The functions $f_i$, $i\in[N]$, are $p$th-order continuously differentiable in $\R^{m\times n}$ for some $p\ge2$, and moreover, there exists some $L_{p,*}>0$ such that 
\begin{align*}
\|\mathcal{D}^p f_i(X) - \mathcal{D}^p f_i(Y)\|_*\le L_{p,*}\|X-Y\|\qquad\forall X,Y\in\R^{m\times n},i\in[N].    
\end{align*}
\end{assumption}

Under Assumption \ref{asp:high-order}, {\it DeMuon-A} sets the extrapolation count $q$ to $p-1$, adapting to the smoothness level of $f$. In addition, we adopt the following notation for use in the subsequent analysis:
\begin{align}\label{def:Lplambda-sum-theta}
    L_{p,\lambda} := \frac{3^{p-1}N^pL_{p,*}}{p!}\bigg[2\bigg(\frac{\sqrt{N}\lambda}{1-\lambda}\bigg)^p + 1\bigg],\quad \theta_{[p]} := \sum_{s=1}^{p-1}\theta_s,\quad \theta_{[p]}^\prime := \sum_{s=1}^{p-1} |\theta_s|.
\end{align}
where $\lambda$ is defined in \eqref{upbd:eig-W}, and $\{\theta_s\}_{s\in[p-1]}$ are inputs of {\it DeMuon-A}.

The following theorem gives an upper bound on the consensus error of $\{\mathbf{X}^{(k)}\}$ generated by {\it DeMuon-A}. Its proof is identical to that of Theorem~\ref{thm:cs-error} and is therefore omitted.

\begin{theorem}[{{\bf consensus error}}]\label{thm:cs-error+}
	Suppose that Assumption \ref{asp:basic} holds. Let $\{\mathbf{X}^{(k)}\}$ be generated by Algorithm \ref{alg:r-msgn-1+} with step size $\eta>0$, and let $\lambda$ be defined in \eqref{upbd:eig-W}. Then it holds that for all $k\ge0$,
	\begin{align}\label{ineq:cs-error+}
		\|\mathbf{X}^{(k)} - \mathbf{1}\otimes\bar{X}^{(k)}\| \le \frac{\sqrt{N}\lambda\eta}{1-\lambda}.
	\end{align}
\end{theorem}

As in \cite{he2025heavytail}, to achieve acceleration by leveraging the higher-order smoothness of $f$, the extrapolation parameters $\{\gamma_s\}$ and the momentum parameters $\{\theta_s\}$ in Algorithm \ref{alg:r-msgn-1+} must satisfy the following conditions:
\begin{align}\label{linear-cond}
& \begin{bmatrix}
    1/\gamma_1 & 1/\gamma_2 & \cdots & 1/\gamma_q \\
    1/\gamma_1^2 & 1/\gamma_2^2 & \cdots & 1/\gamma_q^2 \\
    \vdots & \vdots & \ddots & \vdots \\
    1/\gamma_1^q & 1/\gamma_2^q & \cdots & 1/\gamma_q^q
\end{bmatrix} \begin{bmatrix}
    \theta_1 \\
    \theta_2 \\
    \vdots \\
    \theta_q
\end{bmatrix} = \begin{bmatrix}
    1 \\
    1 \\
    \vdots \\
    1
\end{bmatrix},\quad \sum_{s=1}^q\theta_s \in(0,1),
\end{align}
where the extrapolation count $q=p-1$, with $p$ being the degree of smoothness in Assumption \ref{asp:high-order}.

We next derive an upper bound on the stationarity error for $\{\rmX^{(k)}\}$ generated by {\it DeMuon-A}, whose proof is deferred to Section \ref{sec:pf-thm2+}.

\begin{theorem}[{{\bf stationarity error}}]\label{thm:stat-error+}
	Suppose that Assumption \ref{asp:basic} and \ref{asp:high-order} hold. Let $V$ be given in Assumption~\ref{asp:basic}, and let $\lambda$, $(\Delta_f,L_\lambda)$, and $L_{p,\lambda}$ be defined in \eqref{upbd:eig-W}, \eqref{def:Llambda-objgap}, and \eqref{def:Lplambda-sum-theta}, respectively. Let $\{\rmX^{(k)}\}$ be generated by Algorithm \ref{alg:r-msgn-1+} with inputs $q=p-1$, $\eta$, and $\{(\gamma_s,\theta_s)\}$ satisfying \eqref{linear-cond}. Let $\theta_{[p]}$ and $\theta_{[p]}^\prime$ be defined in \eqref{def:Lplambda-sum-theta}. Then it holds that for all $K\ge1$,
	\begin{align}
		\frac{1}{K} \sum_{k=0}^{K-1} \E[\|\bar{g}(\rmX^{(k)})\|_*] &\le \frac{\Delta_f}{K\eta} +\frac{4L_{p,\lambda}\eta^p(1+\sum_{s=1}^{p-1}(|\theta_s|/\gamma_s^p))}{\theta_{[p]}(1-\theta_{[p]})(1-\lambda)} + \frac{4\sqrt{N} \theta_{[p]}^\prime\|V\|_*}{\sqrt{\theta_{[p]}}(1-\theta_{[p]})(1-\lambda)} + \frac{L_\lambda\eta}{2} \nonumber\\
        &\quad + \frac{2\theta_{[p]}^\prime\lambda\sqrt{N}\|V\|_*}{(1-\theta_{[p]})\sqrt{1-\lambda}} + \frac{2N\lambda L_\lambda\eta\sum_{s=1}^{p-1}(|\theta_s|/\gamma_s)}{K(1-\lambda)(1-\theta_{[p]})}  + \frac{4\|\nabla\rmF(\rmX^{(0)})\|_*}{K\theta_{[p]}(1-\theta_{[p]})(1-\lambda)}.\label{ineq:stat-error+}
	\end{align}
\end{theorem}

The following theorem provides a convergence guarantee for {\it DeMuon-A} when $\eta$, $\{\gamma_s\}$, and $\{\theta_s\}$ are selected properly. Its proof is relegated to Section \ref{subsec:proof-thm-conv+}.

\begin{theorem}[{{\bf convergence}}]\label{thm:conv+}
	Suppose that Assumption \ref{asp:basic} and \ref{asp:high-order} hold. Let $V$ be given in Assumption \ref{asp:basic}, and let $\lambda$, $(\Delta_f,L_\lambda)$, and $L_{p,\lambda}$ be defined in \eqref{upbd:eig-W}, \eqref{def:Llambda-objgap}, and \eqref{def:Lplambda-sum-theta}, respectively. Let $K$ be the maximum iteration threshold of Algorithm \ref{alg:r-msgn-1+} such that $K\ge\frac{16L_{p,\lambda}^{1/p}(1-\lambda)\Delta_f}{(\sqrt{N}\|V\|_*)^{(p+1)/p}}$, and define
	\begin{align}
		&\hat\gamma = \frac{L_{p,\lambda}^{\frac{2}{3p+1}}}{(\sqrt{N}\|V\|_*)^{\frac{2p+2}{3p+1}}}\cdot\bigg[\frac{(1-\lambda)\Delta_f}{K}\bigg]^{\frac{2p}{3p+1}},\quad \hat\eta = \bigg[\frac{(1-\lambda)\Delta_f}{L_{p,\lambda}K}\bigg]^{\frac{1}{p+1}}\cdot \hat{\gamma}^{\frac{p}{p+1}},\label{def:hat-theta-hat-eta+}\\
        &\hat\gamma_s = \frac{\hat\gamma}{s^2},\quad \hat\theta_s = \frac{\prod_{1\le r\le q, r\neq s}(1-r^2/\hat\gamma)}{(s^2/\hat\gamma)\prod_{1\le r\le q,r\neq s}((s^2-r^2)/\hat\gamma)}.\label{def:hat-theta-hat-eta++}
	\end{align}
	Let $\{\rmX^{(k)}\}$ be generated by Algorithm \ref{alg:r-msgn-1+} with inputs $q=p-1$, $\eta=\hat{\eta}$, and $\{(\gamma_s,\theta_s)\}=\{(\hat\gamma_s,\hat\theta_s)\}$. Then, 
	\begin{align}
		%\frac{1}{K} \sum_{k=0}^{K-1} \E[\|\bar{g}(\rmX^{(k)})\|_*] \le \mathcal{O}\bigg(\frac{N^{\frac{5}{8}}\|V\|_*^{\frac{1}{2}}}{1-\lambda}\bigg(\frac{L_*}{K}\bigg)^{\frac{1}{4}}\bigg),\quad \|\rmX^{(k)} - \mathbf{1}\otimes\bar{X}^{(k)}\| \le \mathcal{O}\bigg(\frac{1}{N^{\frac{1}{8}}\|V\|_*^{\frac{1}{2}}L_*^{\frac{1}{4}}K^{\frac{3}{4}}}\bigg). 
        \frac{1}{K} \sum_{k=0}^{K-1} \E[\|\bar{g}(\rmX^{(k)})\|_*] \le \mathcal{O}(K^{-p/(3p+1)}),\quad \max_{0\le k\le K-1}\{\|\rmX^{(k)} - \mathbf{1}\otimes\bar{X}^{(k)}\|\} \le \mathcal{O}(K^{-(2p+1)/(3p+1)}).\label{upbd:cs-ave+}
	\end{align}
\end{theorem}

\begin{remark}
	The dependence on $K$ in the stationarity error in \eqref{upbd:cs-ave+} matches the results in the centralized case as in \cite[Theorem 4]{he2025heavytail} with tail index to be $2$.
\end{remark}

\section{Numerical results}
\label{sec:sr}

In this section, we compare {\it DeMuon} with baselines for %decentralized matrix {\color{blue}regression} and 
decentralized transformer training.

%\subsection{Decentralized matrix optimization {\color{blue}(give the specific name of the problem)}}
%\textcolor{blue}{For decentralized matrix optimization setting, ....}

\begin{figure*}[ht]
  \centering

  %==================== Row 1: Training ====================%
  \subfloat[Complete Graph]{%
    \begin{minipage}[t]{0.32\linewidth}
      \centering
      \includegraphics[width=\linewidth]
      {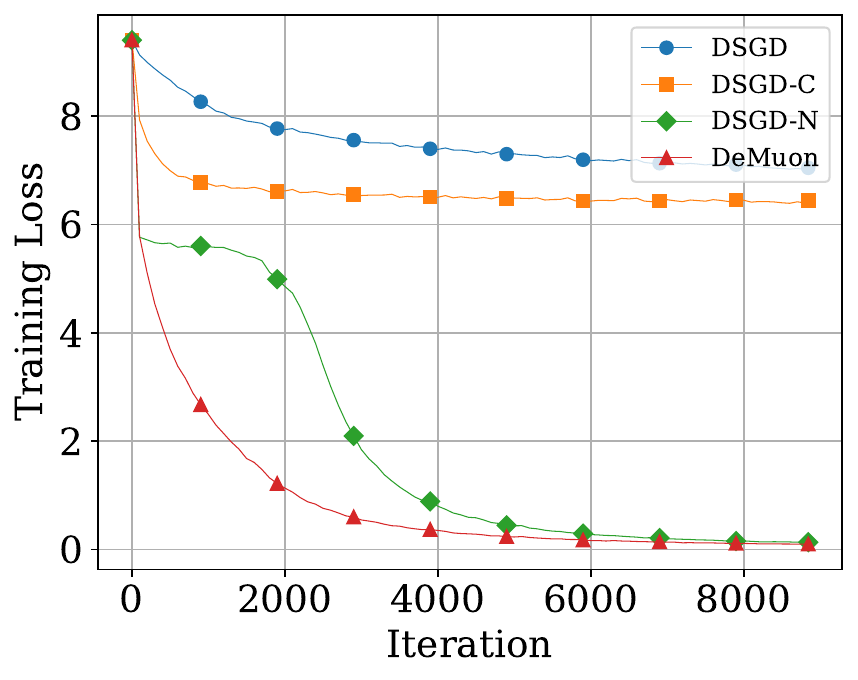}
    \end{minipage}
    \label{fig:complete_train}
  }\hfill
  \subfloat[Directed Exponential Graph]{%
    \begin{minipage}[t]{0.32\linewidth}
      \centering
      \includegraphics[width=\linewidth]
      {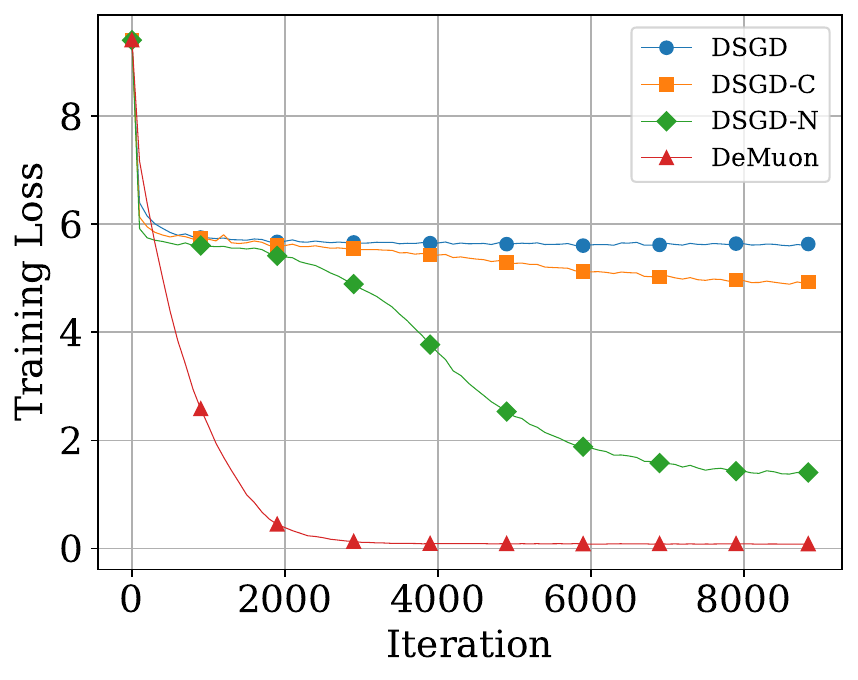}
    \end{minipage}
    \label{fig:exp_train}
  }\hfill
  \subfloat[Ring Graph]{%
    \begin{minipage}[t]{0.32\linewidth}
      \centering
      \includegraphics[width=\linewidth]
      {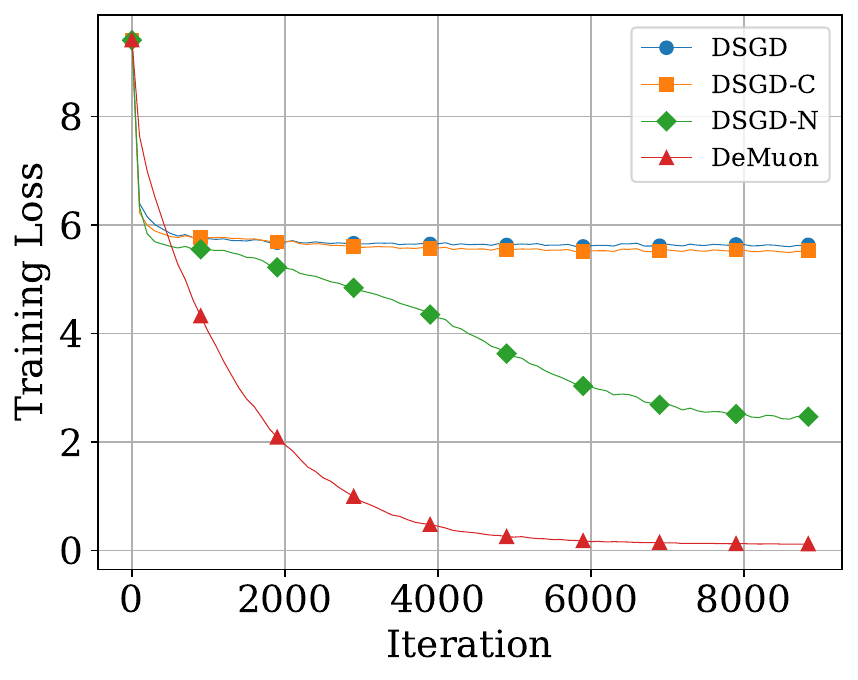}
    \end{minipage}
    \label{fig:ring_train}
  }

  %==================== Row 2: Validation ====================%
  \subfloat[Complete Graph]{%
    \begin{minipage}[t]{0.32\linewidth}
      \centering
      \includegraphics[width=\linewidth]
      {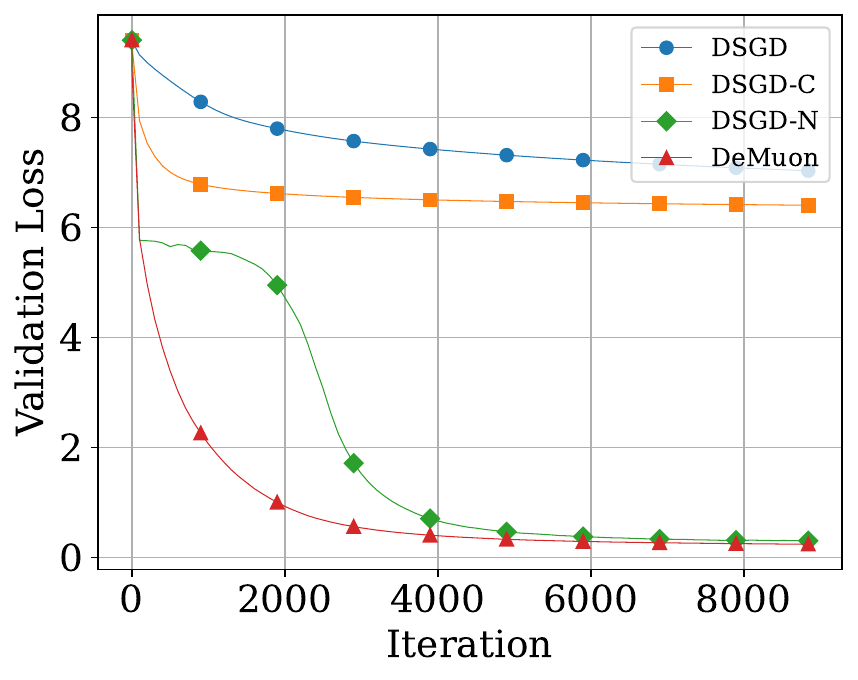}
    \end{minipage}
    \label{fig:complete_val}
  }\hfill
  \subfloat[Directed Exponential Graph]{%
    \begin{minipage}[t]{0.32\linewidth}
      \centering
      \includegraphics[width=\linewidth]
      {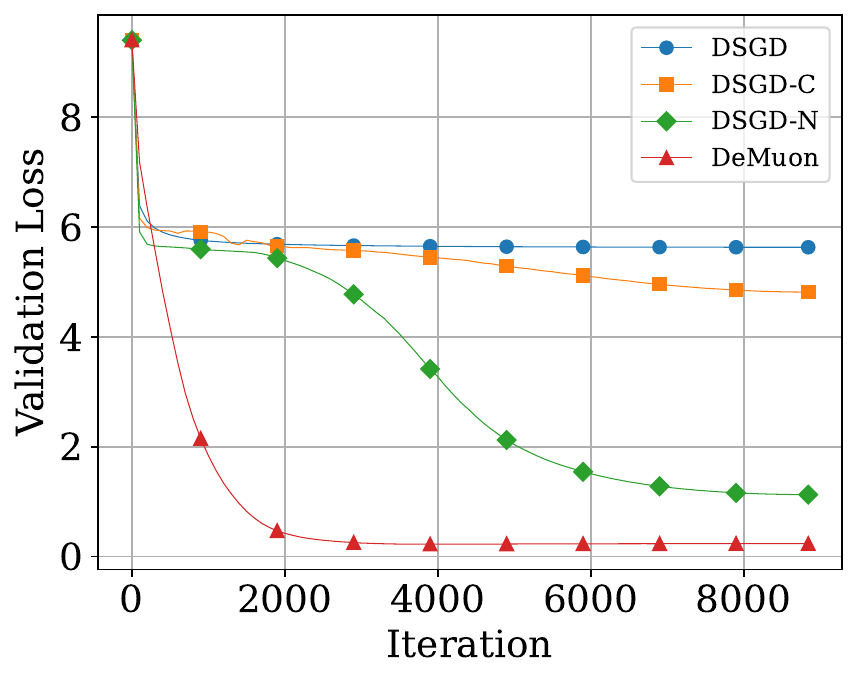}
    \end{minipage}
    \label{fig:exp_val}
  }\hfill
  \subfloat[Ring Graph]{%
    \begin{minipage}[t]{0.32\linewidth}
      \centering
      \includegraphics[width=\linewidth]
      {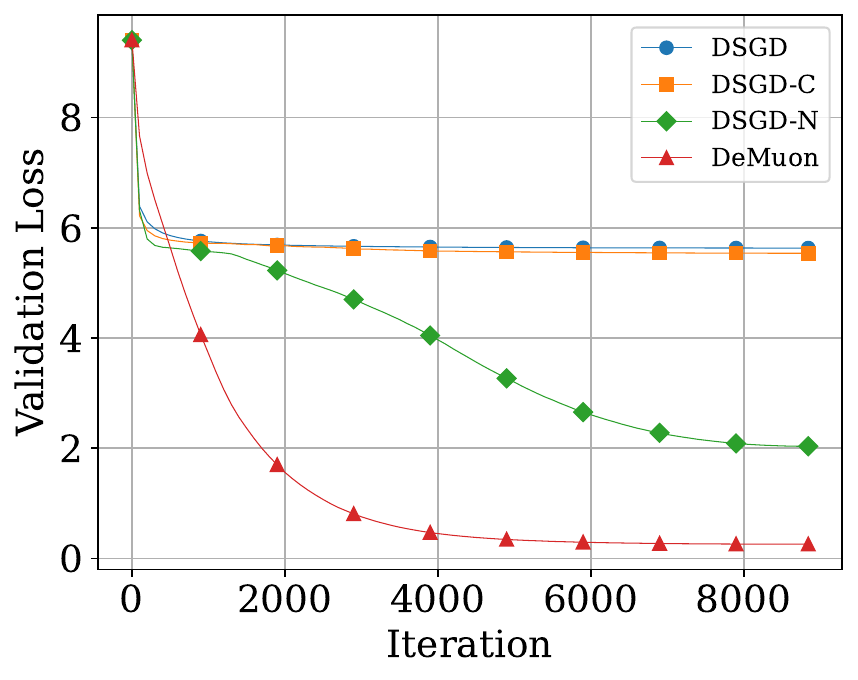}
    \end{minipage}
    \label{fig:ring_val}
  }

  %==================== Row 3: Consensus ====================%
  \subfloat[Complete Graph]{%
    \begin{minipage}[t]{0.32\linewidth}
      \centering
      \includegraphics[width=\linewidth]
      {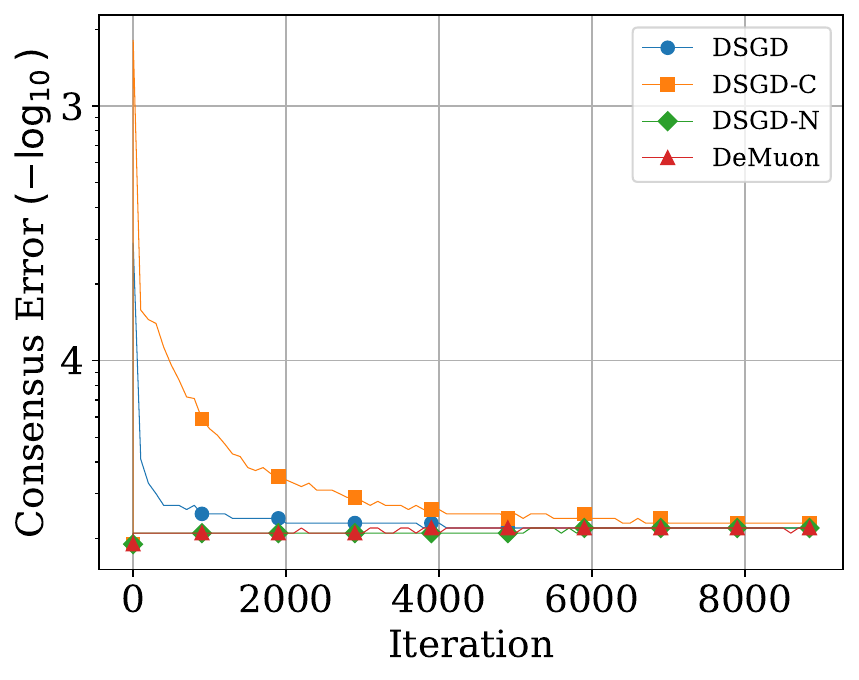}
    \end{minipage}
    \label{fig:complete_cons}
  }\hfill
  \subfloat[Directed Exponential Graph]{%
    \begin{minipage}[t]{0.32\linewidth}
      \centering
      \includegraphics[width=\linewidth]
      {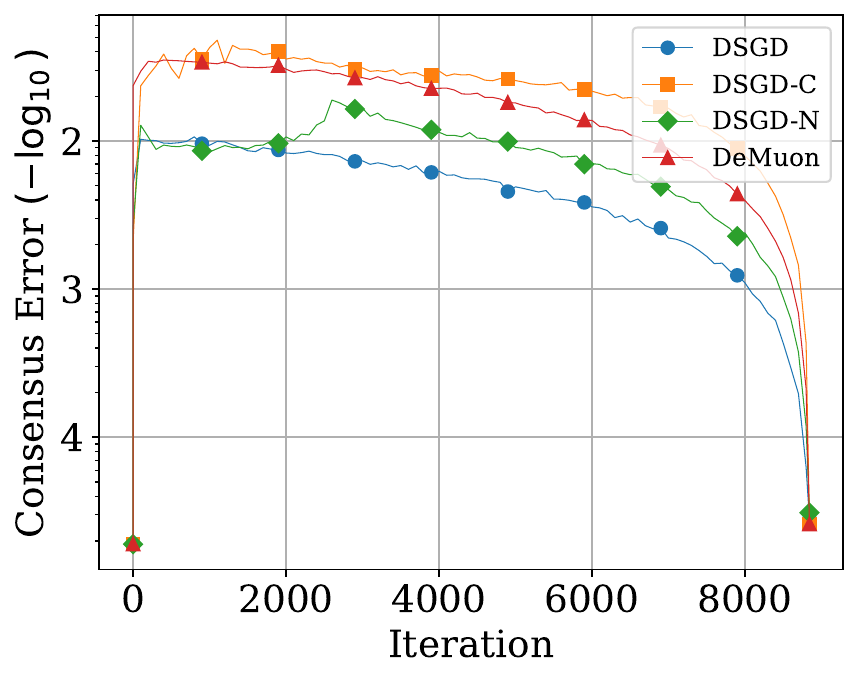}
    \end{minipage}
    \label{fig:exp_cons}
  }\hfill
  \subfloat[Ring Graph]{%
    \begin{minipage}[t]{0.32\linewidth}
      \centering
      \includegraphics[width=\linewidth]
      {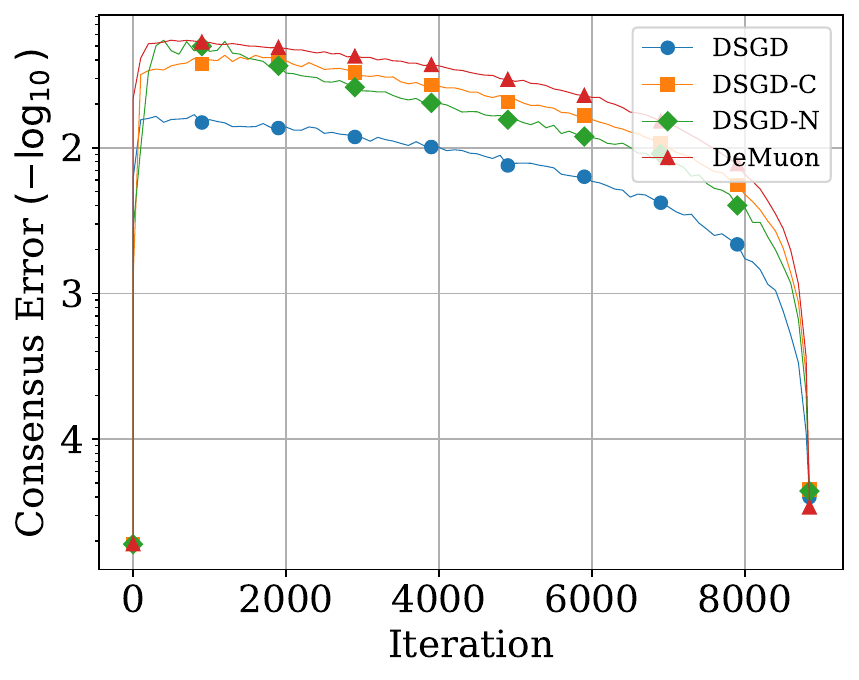}
    \end{minipage}
    \label{fig:ring_cons}
  }

  \caption{
  Training loss (row 1), validation loss (row 2), and consensus error (row 3) in decentralized training of Transformer models over complete (column 1), directed exponential (column 2), and ring (column 3) graph topologies.
  }

  \label{fig:all_results}
\end{figure*}

\subsection{Decentralized transformer training}

In this subsection, we evaluate the performance of {\it DeMuon} on decentralized transformer training. Our results serve to (i) verify that matrix orthogonalization and gradient tracking in {\it DeMuon} translate into practical gains, and (ii) assess robustness to various communication topologies.

%In decentralized deep learning, we evaluate {\it DeMuon} on decentralized training of Transformer-based language models. The experiments serve two purposes: to verify that the matrix orthogonalization and gradient tracking in {\it DeMuon} translate into practical gains, and to assess robustness to the communication topology. %We use a Transformer model because prior work suggests that its training dynamics exhibit heavy-tailed gradient noise, which aligns with the setting of Section~\ref{sec:d-muon}.

{\bf Experiment setup.} We consider auto-regressive language modeling on the English portion of the Multi30k English-German translation dataset \citep{Multi30K}. We adopt the 7.4M-parameter GPT model \citep{radford2018improving} %mainly following \cite{yu2026decentralized} but 
with embedding dimension $256$ and feedforward dimension $1024$. %(versus $128$ and $512$ in \cite{yu2026decentralized}), 
%giving $7.4$M parameters instead of ${\approx}3$M. 
The training dataset split contains $29{,}000$ sentences ($377{,}687$ tokens) and the validation dataset split contains $1{,}014$ sentences ($13{,}326$ tokens); the vocabulary size is $10{,}208$. We use $N=8$ nodes with a batch size of $64$, and training runs for $12$ epochs ($8{,}844$ iterations). 

%The training data (Multi30k train split) is tokenized, shuffled once, and partitioned equally and disjointly across the $N$ nodes so that each node has a fixed subset of the data; validation uses the official dev split jointly. 

We compare {\it DeMuon} with three baseline methods: (i) {\it DSGD}~\citep{nedic2009distributed}, a decentralized vanilla SGD without gradient tracking; (ii) {\it DSGD-C}~\citep{dsgd-clip}, a decentralized clipped SGD without gradient tracking; and (iii) {\it DSGD-N}~\citep{yu2026decentralized}, a decentralized normalized SGD with gradient tracking. 
%Thus, the comparison highlights the benefit of matrix-variate updates via $\mathrm{msgn}(\cdot)$ in {\it DeMuon} (Algorithm~\ref{alg:r-msgn-1}) versus vector-style or non-robust methods.
We also consider three communication graphs: (i) a \emph{complete} graph (full connectivity); (ii) a \emph{directed exponential} graph (one fixed realization per run); and (iii) a \emph{ring} graph (minimal connectivity). These choices of graphs span a range of mixing rates $\lambda$ defined in \eqref{upbd:eig-W}: for the complete graph, we use uniform weights $1/N$, which yields $\lambda = 0$; for the directed exponential and ring graphs, we use Metropolis weights, yielding $\lambda \approx 0.60$ and $\lambda \approx 0.80$, respectively.

For each competing method and communication graph, we tune the hyperparameters so that all methods are compared on a similar consensus-error scale. On the complete graph ($\lambda=0$), {\it DSGD-N} uses a constant learning rate while {\it DSGD}, {\it DSGD-C}, and {\it DeMuon} use diminishing learning rates; on the exponential and ring graphs ($\lambda>0$), all methods use a linear decay schedule, i.e, $\eta(1-k/K)$ at iteration $k$, where $K$ is the maximum iteration threshold. %which maintains a large effective learning rate for most of training and drives the learning rate to zero at the final iteration.
Hyperparameters are tuned per method and topology and are listed in Appendix~\ref{app:params}, Table~\ref{tab:stepsizes}. % For each method on each graph, we report only the result with the lowest consensus error. All runs use random seed 42.

{\bf Main results.} Fig.~\ref{fig:all_results} reports the training and validation losses, and the consensus error over iterations, for each topology. In Table~\ref{tab:results}, we present the validation loss, perplexity, and consensus error at the final iteration for each method and topology. The perplexity is defined as the exponential of the validation loss. The consensus error is computed as the $\ell_2$-norm of the vector whose entries are the spectral norms of the individual weight matrix blocks. 

From Fig.~\ref{fig:all_results} and Table~\ref{tab:results}, we observe that {\it DeMuon} outperforms {\it DSGD} variants across all three communication topologies. This shows that the advantages of {\it Muon} in the centralized setting can be effectively brought into the decentralized setting.

%Table~\ref{tab:results} summarizes the final validation loss, perplexity ($e^{\text{loss}}$), and consensus error for each method and topology.

%Since a neural network consists of multiple parameter matrices of different shapes, one cannot directly form a single stacked matrix and compute its spectral norm as in Theorem~\ref{thm:cs-error}. Instead, we compute the consensus error per parameter: for each parameter, we stack the $N$ workers' deviations from the mean vertically into one matrix and take its spectral norm. The reported consensus error is the $\ell_2$-norm of the vector of per-parameter spectral norms, which upper-bounds the theoretical consensus error in Theorem~\ref{thm:cs-error} by the triangle inequality.

\begin{table}[ht]
\caption{Final validation loss, perplexity, and consensus error after 12 epochs for each method and topology. Lower values are better for all metrics. %Perplexity $= e^{\text{loss}}$ measures the effective number of equally likely next-token choices.
}
\centering
\label{tab:results}
\begin{tabular}{llll}
\hline
{\it Algorithm} & Complete & Directed Exp. & Ring \\
\hline
\multicolumn{4}{l}{{\bf Validation\ Loss\,/\,Perplexity}} \\[2pt]
{\it DSGD}      & 7.031\,/\,1131.24 & 5.632\,/\,279.23 & 5.632\,/\,279.30 \\
{\it DSGD-C} & 6.403\,/\,603.77  & 4.816\,/\,123.44 & 5.538\,/\,254.06 \\
{\it DSGD-N}  & 0.302\,/\,1.35   & 1.127\,/\,3.09   & 2.035\,/\,7.65 \\
{\it DeMuon}    & 0.240\,/\,1.27   & 0.234\,/\,1.26   & 0.258\,/\,1.29 \\
\hline
\multicolumn{4}{l}{\bf Consensus Error} \\[2pt]
{\it DSGD}       & $2.2\!\times\!10^{-5}$ & $2.8\!\times\!10^{-5}$ & $4.0\!\times\!10^{-5}$ \\
{\it DSGD-C} & $2.3\!\times\!10^{-5}$ & $2.6\!\times\!10^{-5}$ & $4.5\!\times\!10^{-5}$ \\
{\it DSGD-N}  & $2.2\!\times\!10^{-5}$ & $3.1\!\times\!10^{-5}$ & $4.4\!\times\!10^{-5}$ \\
{\it DeMuon}     & $2.2\!\times\!10^{-5}$ & $2.6\!\times\!10^{-5}$ & $3.4\!\times\!10^{-5}$ \\
\hline
\end{tabular}
\end{table}

\section{Acknowlegement}

The simulations were enabled by resources provided by the National Academic Infrastructure for Supercomputing in Sweden (NAISS), partially funded by the Swedish Research Council through grant agreement no. 2022-06725. The work of Chuan He was
partially supported by the Wallenberg AI, Autonomous Systems and Software Program (WASP) funded by the Knut and
Alice Wallenberg Foundation. The work of Shuyi Ren and Erik G. Larsson was supported in part by ELLIIT, VR, and the KAW foundation.

\bibliographystyle{tmlr}
\bibliography{ref}

\newpage

\appendix

\section*{Appendix}

\section{Related works}
\label{Appendix:related_works}

We now provide a discussion of related works on matrix-variate optimizers and decentralized optimization.

%%%%%%% Muon
{\bf Muon and variations.} The empirical success of Muon has inspired a large number of research efforts focused on analyzing its convergence \citep{li2025note,shen2025convergence,sato2025analysis,chen2025muon,sfyraki2025lions,kim2026convergence}. Specifically, \cite{li2025note} give one of the first dedicated convergence analyses of Muon and interpret it as spectral-norm steepest descent. \cite{shen2025convergence} provide a convergence analysis of Muon and its comparison with GD, identifying conditions under which Muon outperforms GD. Their results show that Muon can exploit low-rank Hessian structures commonly observed in neural network training. \cite{sato2025analysis} provide theoretical support for Muon by establishing convergence guarantees across practical settings with and without Nesterov momentum and weight decay. Their analysis shows that weight decay leads to tighter theoretical bounds, characterizes its interaction with the learning rate, and further derives the critical batch size that minimizes training cost. \cite{chen2025muon} provide a theoretical understanding of Muon by placing it within the Lion-($\mathcal K$) optimizer family, showing that Muon corresponds to Lion-($\mathcal K$) with the nuclear norm. Their analysis reveals that Muon implicitly enforces spectral norm constraints on weight matrices and further motivates broader implicitly regularized optimization methods through alternative choices of convex map $\mathcal K$. \cite{sfyraki2025lions} unify Lion and Muon under the Stochastic Frank--Wolfe framework, showing that both methods with weight decay admit convergence guarantees via the Frank--Wolfe gap and converge to KKT points under norm constraints. They further develop robust variants for heavy-tailed noise, leading to improved Lion and Muon algorithms with stronger theoretical guarantees and enhanced practical applicability. In recent works, \cite{kim2026convergence} give the first nonconvex convergence analysis for the practical finite-step Newton–Schulz version of Muon. 

%\cite{shah2025practical} argues Muon improves the compute–time Pareto frontier and remains data-efficient at large batch sizes. 

Besides PolarGrad~\citep{lau2025polargrad} and LR-Muon~\citep{he2025low} discussed in Section~\ref{sec:intro}, several other Muon variants have also been developed, including Scion~\citep{pethick2025training}, Gluon~\citep{riabinin2025gluon}, Dion~\citep{ahn2025dion}, EF21-Muon~\citep{gruntkowska2025error}, AdaMuon~\citep{si2025adamuon}, MuLoCo~\citep{therien2025muloco}, MuonEq~\citep{chang2026muoneq}, and Muon+~\citep{zhang2026muon+}, among others. Below, we briefly introduce several representative variants. \cite{pethick2025training} place Muon-like updates in a broader LMO family over norm balls and shows these methods can be used for unconstrained deep learning. \cite{riabinin2025gluon} generalize Muon/Scion under a refined layer-wise smoothness model and narrows the gap between theory and practice for LMO-based optimizers. \cite{ahn2025dion} propose a scalable and efficient update rule that replaces Newton--Schulz iteration with amortized power iteration on a momentum buffer, eliminating full-matrix reconstruction and naturally supporting weight sharding. Their rank-fraction mechanism with error feedback further enables low-rank updates, balancing optimization quality and computational efficiency. \cite{gruntkowska2025error} introduce EF21-Muon, the first communication-efficient error-feedback framework for non-Euclidean LMO-based optimizers. However, all of previous mentioned works are not considering decentralized peer-to-peer setting. \cite{si2025adamuon} combine element-wise second-moment adaptivity with orthogonalized updates and reports more than 40\% training-efficiency gains over Adam in large-scale settings. \cite{therien2025muloco} use Muon as the local inner optimizer in a low-communication distributed pretraining loop and shows strong compressibility with error feedback. \cite{chang2026muoneq} add light row/column equilibration before finite-step orthogonalization and preserves the standard Muon-type stationarity guarantee while improving conditioning. Recently, \cite{zhang2026muon+} add a post-orthogonalization normalization, reporting consistent gains over Muon across GPT- and LLaMA-style pretraining runs. Moreover, \cite{du2026newton} derive a right-preconditioned Muon-like update from a quadratic surrogate and interprets vanilla Muon as an implicit Newton-type method missing right preconditioning. \cite{boissin2025turbo} speed up Newton–Schulz orthogonalization through preconditioning and reports up to 2.8×speedup in the approximation routine with end-to-end training gains. 

%%%%%%% Decentralized
{\bf Decentralized optimization.} Decentralized optimization has been extensively studied over the past decades. Early works focused on fundamental algorithmic frameworks, including distributed subgradient method \citep{nedic2009distributed} and dual averaging \citep{duchi2011dual}, which establish the basis for consensus-based optimization over networks. Building on these foundations, a number of refined first-order methods have been proposed to improve convergence and scalability, such as NEXT \citep{di2016next}, EXTRA \citep{shi2015extra}, and CHOCO-SGD \citep{koloskova2019decentralized}. Although these methods were originally designed for general distributed optimization tasks, they have recently gained renewed attention in the context of large-scale machine learning, where decentralized architectures provide a promising alternative to centralized training for deep neural networks. And there are a lot of variants has been explored in \cite{tang2018communication}.

For decentralized optimization in deep models training, the most established decentralized deep-learning baselines are still decentralized SGD variants. Building upon PSGD~\citep{zinkevich2010parallelized}, D-PSGD \citep{lian2017can} is the first influential result to argue that decentralized training can outperform centralized parameter-server-style training in regimes where communication at the busiest node is the main bottleneck, and its experiments show up to an order-of-magnitude speed advantage under low bandwidth or high latency. 
AD-PSGD \citep{lian2018asynchronous} then extends the picture to heterogeneous environments by combining decentralization with asynchrony, proving the same optimal ($O(1/\sqrt{K})$, where $K$ is the number of iterations) rate as standard SGD together with linear speedup in the number of workers.

The decentralized Adam-family \citep{kingma2015adam, zhang2024adam} literature is more delicate.  \cite{nazari2022dadam} introduce DADAM, a consensus-based distributed adaptive moment method, and provide dynamic-regret guarantees in stochastic and deterministic settings, framing the problem as one of decentralized computation without a central node. But later work on convergent decentralized adaptive methods \citep{chen2023convergence} sharpens the picture substantially: it proposes a general recipe for converting convergent adaptive methods into decentralized counterparts, emphasizes the need for consensus on adaptive learning rates, derives a large-($K$) rate of order ($O(\sqrt{n}/\sqrt{K}$), where $n$ is optimizer's dimension) for decentralized AMSGrad \citep{tran2019convergence}, and explicitly shows that DADAM can diverge because it lacks the right consensus mechanism on adaptive-rate statistics. Shampoo: \cite{gupta2018shampoo} is a structure-aware stochastic preconditioning method for tensor optimization that uses per-dimension matrix preconditioners to achieve faster convergence than standard optimizers while maintaining comparable per-step computational cost. And, its distributed version in \cite{shi2023distributed}, a scalable distributed implementation of Shampoo using block-diagonal Kronecker preconditioning and PyTorch DTensor enables near-parity runtime with standard optimizers while achieving superior training performance on large-scale deep networks.

Moreover, decentralized matrix optimization problem in machine learning is absolutely real—matrix completion, matrix factorization, matrix sensing, and decentralized PCA all have active literature—yet these works are overwhelmingly built around factorized Euclidean/Frobenius or manifold geometries, not operator-norm-aware Muon-style updates. That leaves a natural opening for decentralized spectral-norm-aware matrix optimization. \cite{kempe2004decentralized}'s work is one of the earliest and clearest demonstrations that spectral computation can be done without centralizing the matrix. \cite{gemulla2011large} introduce stratified SGD (SSGD) and specialize it to DSGD for large sparse matrix factorization. \cite{balcan2014improved} use improved sketching/subspace-embedding ideas to reduce communication and computation for distributed PCA. \cite{ling2012decentralized} provide a early decentralized low-rank matrix-completion method using factorization with public/private matrix variables across agents. \cite{mishra2016riemannian} formulate decentralized low-rank matrix completion on the Grassmann manifold with gossip-based consensus. \cite{zhu2019distributed} study DGD+LOCAL for distributed low-rank matrix approximation and shows exact consensus/global-optimality guarantees under favorable geometry. \cite{maros2023decentralized} develop the first decentralized matrix-sensing theory with statistical, communication, and convergence guarantees for a nonconvex Burer–Monteiro approach. \cite{philippenko2025depth} analyze low-rank matrix factorization across clients and it into a strongly convex problem after initialization. \cite{saadati2026decaf} extend decentralized LoRA and introduces TSVD-based factorization to remove consensus interference.

\section{Proof overview}
\label{sec:pf}

In this section, we provide an overview of the proofs in Sections \ref{sec:d-muon} and \ref{sec:DeMuon+}, as well as some technical results.

For convenience, we denote $\xi:=\{\xi_i\}_{i\in[N]}$ and adopt the following notation:
\begin{align}\label{def:not-GOO}
\mathbf{G}(\mathbf{X};\xi) := \left[\begin{matrix}
G_1(X_1;\xi_1)\\
\vdots\\
G_N(X_N;\xi_N)
\end{matrix}\right],\ \
\mathbf{O}:= \left[\begin{matrix}
\mathrm{msgn}(V_1)\\
\vdots\\
\mathrm{msgn}(V_N)
\end{matrix}\right],\ \ \bar{O}:=\frac{1}{N}\sum_{i=1}^N\mathrm{msgn}(V_i).
\end{align}
Using stacked notation, the updates of {\it DeMuon} in \eqref{update-mk}-\eqref{update-xk} can be rewritten compactly as:
\begin{align}\label{def:ave-update}
	\left\{
	\begin{array}{l}
		\rmM^{(k)}  = (1-\theta)\rmM^{(k-1)} + \theta \rmG(\rmX^{(k)};\xi^{(k)}),\\[3pt]
		\rmV^{(k)} = (W\otimes I_m)(\rmV^{(k-1)}+\rmM^{(k)}-\rmM^{(k-1)}),\\[3pt]
		\rmX^{(k+1)} = (W\otimes I_m)(\rmX^{(k)} - \eta \rmO^{(k)}),
	\end{array}
	\right.    \qquad \forall k=0,1,\ldots.
\end{align}
Similarly, the updates of {\it DeMuon-A} in \eqref{update-zk+}-\eqref{update-xk+} can be rewritten as:
\begin{align}\label{def:ave-update+}
	\left\{
	\begin{array}{l}
        \rmZ^{(k,s)}  = \rmX^{(k)} + \frac{1-\gamma_s}{\gamma_s}(\rmX^{(k)} - \rmX^{(k-1)})\quad \forall s\in[p-1],\\[3pt]
		\rmM^{(k)}  = (1-\sum_{s=1}^{p-1}\theta_s)\rmM^{(k-1)} + \sum_{s=1}^{p-1}\theta_s \rmG(\rmZ^{(k,s)};\xi^{(k,s)}),\\[3pt]
		\rmV^{(k)} = (W\otimes I_m)(\rmV^{(k-1)}+\rmM^{(k)}-\rmM^{(k-1)}),\\[3pt]
		\rmX^{(k+1)} = (W\otimes I_m)(\rmX^{(k)} - \eta \rmO^{(k)}),
	\end{array}
	\right.    \qquad \forall k=0,1,\ldots
\end{align}
For convenience, we adopt the notation:
\begin{align}\label{def:DeltaM-DeltaV}
	\Delta \rmM^{(k)} = \nabla \rmF(\rmX^{(k)}) - \rmM^{(k)},\quad \Delta \rmV^{(k)} = \rmV^{(k)} - \mathbf{1}\otimes\bar{V}^{(k)}\qquad\forall k\ge-1,
\end{align}
where $\{(\rmX^{(k)},\rmM^{(k)},\rmV^{(k)})\}$ are generated by Algorithm \ref{alg:r-msgn-1} or \ref{alg:r-msgn-1+}, and we set $\rmX^{(-1)}=\rmX^{(0)}$. By the double stochasticity of $W$ and the update rules of $\{\rmX^{(k)}\}$ in \eqref{def:ave-update} and \eqref{def:ave-update+}, the updates of $\{\bar{X}^{(k)}\}$ can be written as
\begin{align}\label{update-ave-x}
	\bar{X}^{(k+1)} = \bar{X}^{(k)} - \eta \bar{O}^{(k)}\qquad\forall k=0,1,2,\ldots.
\end{align}

Recall the mixed product property of the Kronecker product:
\begin{align}\label{eq:ppt-pd}
	(A\otimes B)(C\otimes D)=(AC)\otimes(BD)  
\end{align}
when $AC$ and $BD$ are defined (i.e., dimensionally compatible). This property will be used in latter analysis.

The following lemma establishes some useful properties of the spectral and nuclear norms.
\begin{lemma}\label{lem:tech-do}
	Suppose that Assumption \ref{asp:basic} holds. Let $\rmX\in\R^{(Nm)\times n}$ be given, and $W\in\R^{N\times N}$ be the mixing matrix. Then, the following statements hold.
	\begin{enumerate}
		\item[{\rm (i)}] $\max_{1\le i \le N}\{\|X_i\|\}\le \|\mathbf{X}\|$ and $\max_{1\le i \le N}\{\|X_i\|_*\}\le \|\rmX\|_*$.
		\item[{\rm (ii)}] $\|\rmX\|\le \sqrt{\sum_{i=1}^N\|X_i\|^2}$ and $\|\rmX\|_*\le\sum_{i=1}^N\|X_i\|_*$
		\item[{\rm (iii)}] $\frac{1}{N}(\mathbf{1}\mathbf{1}^T\otimes I_m)\rmX = \mathbf{1}\otimes\bar{X}$.
		\item[{\rm (iv)}] $(W-\frac{1}{N}\mathbf{1}\mathbf{1}^T)(I_N - \frac{1}{N}\mathbf{1}\mathbf{1}^T)=W-\frac{1}{N}\mathbf{1}\mathbf{1}^T=(I_N - \frac{1}{N}\mathbf{1}\mathbf{1}^T)(W-\frac{1}{N}\mathbf{1}\mathbf{1}^T)$.
	\end{enumerate}
\end{lemma}

\begin{proof}
	To prove statements (i) and (ii), it suffices to show that
	\begin{align}\label{ineq:AB-spec-nuclear}
		\|A\| \le \left\|\left[\begin{matrix} A \\ B\end{matrix}\right]\right\|\le \sqrt{\|A\|^2 + \|B\|^2},\quad \|A\|_*\le\left\|\left[\begin{matrix} A \\ B\end{matrix}\right]\right\|_* \le \|A\|_* + \|B\|_*  
	\end{align}
	hold for all $A\in\R^{m^\prime\times n},B\in\R^{m^{\prime\prime}\times n}$. To this end, we denote $P:=[I_{m^\prime},\mathbf{0}]\in\R^{m^\prime\times(m^\prime+m^{\prime\prime})}$. Then, one has
	%Using the (mixed) submultiplicativity of spectral norm and nuclear norm, we obtain that  
	\begin{align}
		\|A\| = \left\|P\left[\begin{matrix} A \\ B\end{matrix}\right]\right\|\le \|P\|\cdot\left\|\left[\begin{matrix} A \\ B\end{matrix}\right]\right\|=\left\|\left[\begin{matrix} A \\ B\end{matrix}\right]\right\|,\quad \|A\|_* = \left\|P\left[\begin{matrix} A \\ B\end{matrix}\right]\right\|_*\le \|P\|\cdot\left\|\left[\begin{matrix} A \\ B\end{matrix}\right]\right\|_*=\left\|\left[\begin{matrix} A \\ B\end{matrix}\right]\right\|_*.\label{two-sp-un-upbd-1}  
	\end{align}
	where the inequalities are due to the submultiplicativity of the spectral norm and the nuclear norm, respectively. In addition, notice that
	\begin{align*}
		\left\|\left[\begin{matrix} A \\ B\end{matrix}\right]\right\|^2  \le \|A^TA\| + \|B^TB\| = \|A\|^2 + \|B\|^2,\quad 
		\left\|\left[\begin{matrix} A \\ B\end{matrix}\right]\right\|_* \le \left\|\left[\begin{matrix} A \\ \mathbf{0}\end{matrix}\right]\right\|_* + \left\|\left[\begin{matrix} \mathbf{0} \\ B\end{matrix}\right]\right\|_* = \|A\|_* + \|B\|_*,
	\end{align*}
	which together with \eqref{two-sp-un-upbd-1} implies \eqref{ineq:AB-spec-nuclear}. Hence, statements (i) and (ii) hold.
	
	Statement (iii) holds because
	\begin{align*}
		\frac{1}{N}(\mathbf{1}\mathbf{1}^T\otimes I_m)\rmX=\frac{1}{N}\begin{bmatrix}
			I_m&\cdots&I_m\\
			\vdots&\ddots&\vdots\\
			I_m&\cdots&I_m
		\end{bmatrix}\cdot\begin{bmatrix}
			X_1\\
			\vdots\\
			X_N
		\end{bmatrix} = \mathbf{1}\otimes\bar{X}.
	\end{align*}
	In addition, statement (iv) holds because $W\mathbf{1}\mathbf{1}^T=\mathbf{1}\mathbf{1}^T=\mathbf{1}\mathbf{1}^TW$ due to Assumption \ref{asp:basic}(d). 
\end{proof}

The following lemma provides a useful property of nuclear norm, which is adapted from \cite[Lemma 8]{an2025asgo}.
\begin{lemma}\label{lem:tech-nuclear-ineq}
	Let $A\in\R^{m\times n}$ be given. Then,
	\begin{align*}
		\|A\|_* \le \sqrt{\|\Lambda\|_*\mathrm{Tr}(A\Lambda^{-1}A^T)}\qquad\forall \Lambda\succ\mathbf{0}.
	\end{align*}
\end{lemma}

With the above technical lemma, we next provide a matrix martingale moment inequality in the following lemma, which will play a central role in our analysis.

\begin{lemma}\label{lem:concentration}
	Let $\{A_t\}_{t=0}^k$ be a set of $p\times q$ random matrices and $\{\Omega_t\}_{t=0}^k$ be a set of $p\times p$ deterministic symmetric matrices. Assume that $\E[A_t\ |\ \mathcal{F}_{t-1}]=0$ for each $1\le t\le k$ and $\E[A_0]=0$, where $\mathcal{F}_t=\sigma(A_0,\ldots,A_t)$ is the natural filtration. Moreover, assume that there exists a $U\in\R^{p\times q}$ such that $\E[A_t^TA_t]\preceq U^TU$ for each $0\le t\le k$. Then, 
	\begin{align*}
		\E\Big[\Big\|\sum_{t=0}^k\Omega_tA_t\Big\|_*\Big] \le \sqrt{\sum_{t=0}^k\|\Omega_t\|^2} \cdot \|U\|_*.
	\end{align*}
\end{lemma}

\begin{proof}
	By Lemma \ref{lem:tech-nuclear-ineq} with $A=\sum_{t=0}^k\Omega_tA_t$, the following holds for any $\Lambda\succ0$:
	\begin{align*}
		\Big\|\sum_{t=0}^k\Omega_tA_t\Big\|_* & \le \sqrt{\|\Lambda\|_*\mathrm{Tr}\Big(\Big(\sum_{t=0}^k\Omega_tA_t\Big)\Lambda^{-1}\Big(\sum_{t=0}^k\Omega_tA_t\Big)^T\Big)} \\
		& = \sqrt{\|\Lambda\|_*\mathrm{Tr}\Big(\sum_{t=0}^k(\Omega_tA_t)\Lambda^{-1}(\Omega_tA_t)^T + \sum_{s\neq t}(\Omega_tA_t)\Lambda^{-1}(\Omega_sA_s)^T\Big)}.
	\end{align*}
	Taking the expectation on both sides of this inequality, we have that the following holds for any $\Lambda\succ0$:
	\begin{align*}
		\E\Big[\Big\|\sum_{t=0}^k\Omega_tA_t\Big\|_*\Big] & \le \E\Bigg[\sqrt{\|\Lambda\|_*\mathrm{Tr}\Big(\sum_{t=0}^k(\Omega_tA_t)\Lambda^{-1}(\Omega_tA_t)^T + \sum_{s\neq t}(\Omega_tA_t)\Lambda^{-1}(\Omega_sA_s)^T\Big)}\Bigg] \\
		& \le \sqrt{\|\Lambda\|_*\bigg\{\sum_{t=0}^k\E[\mathrm{Tr}((\Omega_tA_t)\Lambda^{-1}(\Omega_tA_t)^T)] + \sum_{s\neq t}\E[\mathrm{Tr}((\Omega_tA_t)\Lambda^{-1}(\Omega_sA_s)^T)]\bigg\}} \\
		& = \sqrt{\|\Lambda\|_*\sum_{t=0}^k\E[\mathrm{Tr}(\Omega_t^2A_t\Lambda^{-1}A_t^T)]} \le \sqrt{\|\Lambda\|_*\sum_{t=0}^k\|\Omega_t\|^2\E[\mathrm{Tr}(A_t\Lambda^{-1}A_t^T)]} \\
		& = \sqrt{\|\Lambda\|_*\sum_{t=0}^k\|\Omega_t\|^2\mathrm{Tr}(\E[A_t^TA_t]\Lambda^{-1})} \le \sqrt{\|\Lambda\|_*\sum_{t=0}^k\|\Omega_t\|^2\mathrm{Tr}(U^TU\Lambda^{-1})},
	\end{align*}
	where the second inequality follows from Jensen's inequality because $\sqrt{\cdot}$ is concave, the first equality is due to $\E[A_t\ |\mathcal{F}_{t-1}]=0$ for all $1\le t\le k$ and $\E[A_0]=0$, the third inequality follows from trace H\"older inequality $\mathrm{Tr}(\Omega_t^2A_t\Lambda^{-1}A_t^T)\le\|\Omega_t^2\|\cdot\|A_t\Lambda^{-1}A_t^T\|_*=\|\Omega_t\|^2\cdot\mathrm{tr}(A_t\Lambda^{-1}A_t^T)$, and the last inequality follows from $\E[A_t^TA_t]\preceq U^TU$ for each $t$. Taking $\Lambda =(U^TU)^{1/2}$, we obtain that
	\begin{align*}
		\E\Big[\Big\|\sum_{t=0}^k\Omega_tA_t\Big\|_*\Big] \le \sqrt{\sum_{t=0}^k\|\Omega_t\|^2}\cdot \|U\|_*.
	\end{align*}    
	Hence, this lemma holds as desired.
\end{proof}

\section{Proof of the main results in Section \ref{sec:d-muon}}

In this section, we provide the proofs of our main results, particularly, Theorems \ref{thm:cs-error}-\ref{thm:conv}.

\subsection{Proof of Theorem \ref{thm:cs-error}}\label{sec:pf-thm1}
\begin{proof}
	When $k=0$, this theorem holds since $X_i^{(0)}=\bar{X}^{(0)}$ for all $i\in[N]$. We next prove this theorem for any $k\ge1$. It follows from the updates of $\{\rmX^{(k)}\}$ in \eqref{def:ave-update}, and Lemma \ref{lem:tech-do}(iii) and (iv) that
	\begin{align*}
		\mathbf{X}^{(k)} - \mathbf{1}\otimes \bar{X}^{(k)} &   = \Big(\Big(I_N- \frac{1}{N}\mathbf{1}\mathbf{1}^T\Big)\otimes I_m\Big) \mathbf{X}^{(k)}\\
		&\overset{\eqref{def:ave-update}}{=}\Big(\Big(I_N- \frac{1}{N}\mathbf{1}\mathbf{1}^T\Big)\otimes I_m\Big)(W\otimes I_m)\big(\mathbf{X}^{(k-1)} - \eta \mathbf{O}^{(k-1)}\big) \\
		& = \Big(\Big(W- \frac{1}{N}\mathbf{1}\mathbf{1}^T\Big)\otimes I_m\Big)\big(\mathbf{X}^{(k-1)} - \eta \mathbf{O}^{(k-1)}\big)\\
		& = \Big(\Big(W- \frac{1}{N}\mathbf{1}\mathbf{1}^T\Big)\otimes I_m\Big)\Big(\Big(\Big(I_{N} - \frac{1}{N}\mathbf{1}\mathbf{1}^T\Big)\otimes I_m\Big)\mathbf{X}^{(k-1)} - \eta \mathbf{O}^{(k-1)}\Big)\\
		&=\Big(\Big(W- \frac{1}{N}\mathbf{1}\mathbf{1}^T\Big)\otimes I_m\Big)\big(\mathbf{X}^{(k-1)} - \mathbf{1}\otimes\bar{X}^{(k-1)} - \eta \mathbf{O}^{(k-1)}\big),
	\end{align*}
	where the first equality follows from Lemma \ref{lem:tech-do}(iii) and $I_N\otimes I_m = I_{Nm}$, the third equality follows from \eqref{eq:ppt-pd} and Assumption \ref{asp:basic}(d), the fourth equality is due to \eqref{eq:ppt-pd} and Lemma \ref{lem:tech-do}(iv), and the last equality follows from Lemma \ref{lem:tech-do}(iii). In addition, notice from \eqref{def:not-GOO} and Lemma \ref{lem:tech-do}(ii) that $\|\mathbf{O}^{(k-1)}\|\le(\sum_{i=1}^N\|\mathrm{msgn}(V_i^{(k-1)})\|^2)^{1/2}=\sqrt{N}$. By this, \eqref{upbd:eig-W}, and the above inequality, one has the following holds for all $k\ge1$,
	\begin{align*}
		\|\mathbf{X}^{(k)} - \mathbf{1}\otimes \bar{X}^{(k)}\| & \le \Big\|\Big(W - \frac{1}{N}\mathbf{1}\mathbf{1}^T\Big)\otimes I_m\Big\| \cdot  \big(\|\mathbf{X}^{(k-1)} - \mathbf{1}\otimes \bar{X}^{(k-1)}\| + \eta\|\mathbf{O}^{(k-1)}\|\big)\\
		&\le \Big\|W - \frac{1}{N}\mathbf{1}\mathbf{1}^T\Big\| \cdot  \big(\|\mathbf{X}^{(k-1)} - \mathbf{1}\otimes \bar{X}^{(k-1)}\| + \sqrt{N}\eta\big)\\
		&\overset{\eqref{upbd:eig-W}}{=}\lambda\big(\|\mathbf{X}^{(k-1)} - \mathbf{1}\otimes \bar{X}^{(k-1)}\| + \sqrt{N}\eta\big)
		%\lambda \big(\|X_k - \mathbf{1}\otimes \bar{X}_k\| + \sqrt{N}\eta\big).
	\end{align*}
	where the second inequality is due to $\|A\otimes B\|\le \|A\|\cdot\|B\|$ for all $A\in\R^{m\times n}$ and $B\in\R^{p\times q}$. This recursion implies that
	\begin{align*}
		\|\mathbf{X}^{(k)} - \mathbf{1}\otimes \bar{X}^{(k)}\| \le \lambda^k\|\mathbf{X}^{(0)} - \mathbf{1}\otimes \bar{X}^{(0)}\| + \sqrt{N}\eta \sum_{t=1}^k\lambda^t \le \frac{\sqrt{N}\lambda\eta}{1-\lambda}    
	\end{align*}
	where the last inequality is due to $X^{(0)}_i=\bar{X}^{(0)}$ for each $i$ and $\sum_{t=1}^k\lambda^t\le\frac{\lambda}{1-\lambda}$. Hence, the conclusion of this theorem holds. 
\end{proof}

\subsection{Proof of Theorem \ref{thm:stat-error}}\label{sec:pf-thm2}

The following lemma provides a descent inequality on the network average $\{\Bar{X}^{(k)}\}$ generated by Algorithm \ref{alg:r-msgn-1}.
\begin{lemma}\label{lem:desc-ave}
	Suppose that Assumption \ref{asp:basic} holds. Let $\{(\rmX^{(k)},\rmM^{(k)},\rmV^{(k)})\}$ be the sequence generated by Algorithm \ref{alg:r-msgn-1} with inputs $\eta$ and $\theta$, and let $L_\lambda $ and $\{(\Delta \rmM^{(k)},\Delta \rmV^{(k)})\}$ be defined in \eqref{def:Llambda-objgap} and \eqref{def:DeltaM-DeltaV}, respectively. Then,
	\begin{align}
		f(\bar{X}^{(k+1)}) \le f(\bar{X}^{(k)}) - \eta\|\bar{g} ({\rmX}^{(k)})\|_* + 2\eta(\|\Delta \rmM^{(k)}\|_* + \|\Delta \rmV^{(k)}\|_*)   + \frac{L_\lambda \eta^2}{2}\qquad\forall k\ge0.\label{ineq:desc-1}
	\end{align}
\end{lemma}

\begin{proof}
	Using the update of $\{\rmV^{(k)}\}$ in \eqref{def:ave-update} and Assumption \ref{asp:basic}(d), we have that $\bar{V}^{(k)}=\bar{V}^{(k-1)}+\Bar{M}^{(k)} - \Bar{M}^{(k-1)}$ holds for all $k\ge0$. Observe from Algorithm \ref{alg:r-msgn-1} that $\bar{V}^{(-1)}=\bar{M}^{(-1)}=\mathbf{0}$. In view of these, one can show by induction that $\bar{M}^{(k)}=\bar{V}^{(k)}$ for all $k\ge0$. We next prove \eqref{ineq:desc-1} for any fixed $k\ge0$. By Assumption~\ref{asp:basic}(b) and Lemma \ref{lem:tech-do}(ii), one has  
	\begin{align}
		\|\nabla f(\bar{X}^{(k)}) - \bar{g}({\rmX}^{(k)})\|_* & \le\frac{1}{N}\sum_{i=1}^N\|f_i(\bar{X}^{(k)}) - \nabla f_i(X_i^{(k)})\|_* \nonumber\\
        & \le \frac{L_*}{N}\sum_{i=1}^N\|\bar{X}^{(k)} - X_i^{(k)}\| \le L_*\|\rmX^{(k)} - \mathbf{1}\otimes\bar{X}^{(k)}\|,\label{ineq:ave-grad-1}
	\end{align}
	where the first inequality is due to triangular inequality, the second inequality follows from Assumption~\ref{asp:basic}(b), and the last inequality is because of Lemma \ref{lem:tech-do}(i). Observe that $\|\bar{O}^{(k)}\|\le\frac{1}{N}\sum_{i=1}^N\|\mathrm{msgn}(V^{(k)}_i)\|\le 1$. Using these, \eqref{ineq:desc} with $(X,Y)=(\bar{X}^{(k)},\bar{X}^{(k+1)})$, and \eqref{update-ave-x}, we have
	\begin{align}
		f(\bar{X}^{(k+1)}) & \overset{\eqref{ineq:desc}}{\le} f(\bar{X}^{(k)}) + \langle\nabla f(\bar{X}^{(k)}), \bar{X}^{(k+1)}-\bar{X}^{(k)}\rangle + \frac{L_*}{2}\|\bar{X}^{(k+1)}-\bar{X}^{(k)}\|^2\nonumber\\
		&\overset{\eqref{update-ave-x}}{=} f(\bar{X}^{(k)}) - \eta\langle \bar{V}^{(k)}, \bar{O}^{(k)}\rangle+ \eta\langle \bar{V}^{(k)} - \nabla f(\bar{X}^{(k)}), \bar{O}^{(k)}\rangle + \frac{L_*\eta^2}{2}\|\bar{O}^{(k)}\|^2\nonumber\\
		& \le  f(\bar{X}^{(k)}) - \eta\langle \bar{V}^{(k)}, \bar{O}^{(k)}\rangle + \eta\|\nabla f(\bar{X}^{(k)}) - \bar{V}^{(k)}\|_* + \frac{L_*\eta^2}{2}\nonumber\\
		&\le f(\bar{X}^{(k)}) - \eta\langle \bar{V}^{(k)}, \bar{O}^{(k)}\rangle + \eta\|\nabla f(\bar{X}^{(k)}) - \bar{g}({\rmX}^{(k)})\|_*  + \eta\|\bar{g}({\rmX}^{(k)}) - \bar{V}^{(k)}\|_* + \frac{L_*\eta^2}{2}\nonumber\\
		&\overset{\eqref{ineq:ave-grad-1}}{\le} f(\bar{X}^{(k)}) - \eta\langle \bar{V}^{(k)}, \bar{O}^{(k)}\rangle + L_*\eta\|\rmX^{(k)} - \mathbf{1}\otimes\bar{X}^{(k)}\|  + \eta\|\bar{g}({\rmX}^{(k)}) - \bar{M}^{(k)}\|_* + \frac{L_*\eta^2}{2}\nonumber\\
		&\le f(\bar{X}^{(k)}) - \eta\langle \bar{V}^{(k)}, \bar{O}^{(k)}\rangle  + \eta\|\bar{g}({\rmX}^{(k)}) - \bar{M}_k\|_* + \frac{L_\lambda \eta^2}{2},\label{ineq:upbd-desc-1}
	\end{align}
	where the second inequality follows from $\|\bar{O}^{(k)}\|\le1$ and the trace Hölder inequality, the third inequality is due to the triangular inequality, the fourth inequality follows from \eqref{ineq:ave-grad-1} and $\bar{M}^{(k)}=\bar{V}^{(k)}$, and the last inequality follows from Theorem \ref{thm:cs-error} and \eqref{def:Llambda-objgap}. Note that 
	\begin{align*}
		-\langle\bar{V}^{(k)}, \bar{O}^{(k)}\rangle & = -\frac{1}{N}\sum_{i=1}^N\langle \bar{V}^{(k)}-V^{(k)}_i, \mathrm{msgn}(V^{(k)}_i)\rangle - \frac{1}{N}\sum_{i=1}^N\|V^{(k)}_i\|_*\nonumber\\
		&\le \frac{1}{N}\sum_{i=1}^N\|\bar{V}^{(k)} - V^{(k)}_i\|_* - \frac{1}{N}\sum_{i=1}^N\|V^{(k)}_i\|_* \le \frac{2}{N}\sum_{i=1}^N\|\bar{V}^{(k)} - V^{(k)}_i\|_* - \|\bar{V}^{(k)}\|_*\nonumber\\
		& \le -\|\bar{g}({\rmX}^{(k)})\|_* + \|\bar{g}({\rmX}^{(k)}) - \bar{V}^{(k)}\|_* + \frac{2}{N}\sum_{i=1}^N\|\bar{V}^{(k)} - V^{(k)}_i\|_*,
		%&\overset{\eqref{def:sa-not-2}\eqref{def:sa-not-1}}{\le} -\|\bar{\nabla} F({X}^k)\|_* + \frac{1}{N}\sum_{i=1}^N\|\nabla f_i({X}_i^k) - V_i^k\|_* + \frac{2}{N}\sum_{i=1}^N\|\bar{V}^k - V_i^k\|_*,
	\end{align*}
	where the first inequality follows from the trace H\"older inequality and $\|\mathrm{msgn}(V^{(k)}_i)\|\le1$ for each $i$, and the last two inequalities follow from the triangular inequality. Using this inequality, \eqref{ineq:upbd-desc-1}, Lemma \ref{lem:tech-do}(i), and $\bar{M}^{(k)}=\bar{V}^{(k)}$, we obtain that
	\begin{align}
		f(\bar{X}^{(k+1)}) &  \le f(\bar{X}^{(k)}) - \eta\|\bar{g}({\rmX}^{(k)})\|_* + 2\eta\|\bar{g}({\rmX}^{(k)}) - \bar{M}^{(k)}\|_* + \frac{2\eta}{N}\sum_{i=1}^N\|\bar{V}^{(k)} - V^{(k)}_i\|_* + \frac{L_\lambda \eta^2}{2}\nonumber \\
		&\le f(\bar{X}^{(k)}) - \eta\|\bar{g}({\rmX}^{(k)})\|_* + \frac{2\eta}{N}\sum_{i=1}^N\|\nabla f_i(X_i^{(k)}) - M_i^{(k)}\|_* + \frac{2\eta}{N}\sum_{i=1}^N\|\bar{V}^{(k)} - V^{(k)}_i\|_* + \frac{L_\lambda \eta^2}{2}\nonumber \\
		&\le  f(\bar{X}^{(k)}) - \eta\|\bar{g}({\rmX}^{(k)})\|_* + 2\eta \|\nabla \rmF(\rmX^{(k)}) - \rmM^{(k)}\|_* + 2\eta\|\rmV^{(k)} - \mathbf{1}\otimes \bar{V}^{(k)}\|_*  + \frac{L_\lambda \eta^2}{2}\nonumber\\
		& \overset{\eqref{def:DeltaM-DeltaV}}{=} f(\bar{X}^{(k)}) - \eta\|\bar{g}({\rmX}^{(k)})\|_* + 2\eta\|\Delta \rmM^{(k)}\|_* +2\eta\|\Delta \rmV^{(k)}\|_* + \frac{L_\lambda \eta^2}{2},\nonumber%\label{ineq:desc-interm}
	\end{align}
	where the second inequality is due to the triangular inequality and the third inequality is due to Lemma \ref{lem:tech-do}(i). Hence, this lemma holds as desired.
\end{proof}

The following lemma provides an upper bound on the expected consensus error for the sequence $\{\rmV^{(k)}\}$ generated by Algorithm \ref{alg:r-msgn-1}.
\begin{lemma}
	Suppose that Assumption \ref{asp:basic} holds. Let $\{(\rmM^{(k)},\rmV^{(k)})\}$ be generated by Algorithm \ref{alg:r-msgn-1} with input parameters $\eta$ and $\theta$, let $V$ is given in Assumption \ref{asp:basic}(c), and let $\lambda$ and $\{(\Delta \rmM^{(k)},\Delta \rmV^{(k)})\}$ are defined as in \eqref{upbd:eig-W} and \eqref{def:DeltaM-DeltaV}, respectively. Then it holds that for all $K\ge1$,
    \begin{align}\label{ineq:expec-DVk-bd}
		\sum_{k=0}^{K-1}\E[\|\Delta \rmV^{(k)}\|_*] \le \frac{K\theta\lambda\sqrt{N}\|V\|_*}{(1-\theta)\sqrt{1-\lambda}} + \frac{\theta\lambda\sum_{k=0}^{K-1}\E[\|\Delta \rmM^{(k)}\|_*]}{(1-\theta)(1-\lambda)}.
	\end{align}
\end{lemma}

\begin{proof}
	Fix any $k\ge0$. By \eqref{def:ave-update}, \eqref{def:DeltaM-DeltaV}, \eqref{eq:ppt-pd}, and Lemmas \ref{lem:tech-do}(iii) and (iv), one has that
	\begin{align}
		\Delta \rmV^{(k)} &  \overset{\eqref{def:DeltaM-DeltaV}}{=} \rmV^{(k)} - \mathbf{1}\otimes\bar{V}^{(k)} = \Big(\Big(I_N - \frac{1}{N}\mathbf{1}\mathbf{1}^T\Big)\otimes I_m\Big)\rmV^{(k)}\nonumber\\
		&\overset{\eqref{def:ave-update}}{=}\Big(\Big(I_N - \frac{1}{N}\mathbf{1}\mathbf{1}^T\Big)\otimes I_m\Big)(W\otimes I_m)(\rmV^{(k-1)} +  \rmM^{(k)} - \rmM^{(k-1)}) \nonumber\\
		& \overset{\eqref{eq:ppt-pd}}{=} \Big(\Big(W - \frac{1}{N}\mathbf{1}\mathbf{1}^T\Big)\otimes I_m\Big)(\rmV^{(k-1)} +  \rmM^{(k)} - \rmM^{(k-1)})\nonumber\\
		& \overset{\eqref{def:DeltaM-DeltaV}}{=} \Big(\Big(W - \frac{1}{N}\mathbf{1}\mathbf{1}^T\Big)\otimes I_m\Big)\big(\Delta \rmV^{(k-1)} + \rmM^{(k)} - \rmM^{(k-1)}\big)\nonumber \\
		& = \Big(\Big(W - \frac{1}{N}\mathbf{1}\mathbf{1}^T\Big)\otimes I_m\Big)^{k+1}\Delta \rmV^{(-1)} + \sum_{t=0}^k\Big(\Big(W - \frac{1}{N}\mathbf{1}\mathbf{1}^T\Big)\otimes I_m\Big)^{k-t+1}\big(\rmM^{(t)} - \rmM^{(t-1)}\big)\nonumber\\
		& = \sum_{t=0}^k\Big(\Big(W - \frac{1}{N}\mathbf{1}\mathbf{1}^T\Big)\otimes I_m\Big)^{k-t+1}\big(\rmM^{(t)} - \rmM^{(t-1)}\big),\label{ineq:V-cs-1}
	\end{align}
	where the second equality follows from Lemma \ref{lem:tech-do}(iii), the fourth equality follows from \eqref{eq:ppt-pd} and Assumption~\ref{asp:basic}(d), the fifth equality is due to \eqref{def:DeltaM-DeltaV} and Lemma \ref{lem:tech-do}(iii), and the last equality is due to $V^{(-1)}_i=\mathbf{0}$ for each $i\in[N]$. We also recall from the update of $\{\rmM^{(k)}\}$ in \eqref{def:ave-update} that 
	\begin{align*}
		\rmM^{(t)} - \rmM^{(t-1)}
		= \theta (\rmG(\rmX^{(t)};\xi^{(t)}) - \nabla \rmF(\rmX^{(t)}))  + \theta (\rmM^{(t)} - \rmM^{(t-1)}) +  \theta (\nabla \rmF(\rmX^{(t)}) - \rmM^{(t)})\qquad\forall 0\le t\le k,
	\end{align*}
	which, by rearranging terms, is equivalent to
	\begin{align}\label{eq:Mk-Mk-1}
		\rmM^{(t)} - \rmM^{(t-1)} = \frac{\theta }{1-\theta}(\rmG(\rmX^{(t)};\xi^{(t)}) - \nabla \rmF(\rmX^{(t)})) + \frac{\theta }{1-\theta }(\nabla \rmF(\rmX^{(t)}) - \rmM^{(t)})\qquad\forall 0\le t\le k.
	\end{align}
	Taking the nuclear norm on both sides of \eqref{ineq:V-cs-1}, and using \eqref{upbd:eig-W}, \eqref{def:DeltaM-DeltaV}, \eqref{eq:Mk-Mk-1}, Lemma \ref{lem:tech-do}(i), and Assumption~\ref{asp:basic}(c), we derive that 
	\begin{align}
		\|\Delta \rmV^{(k)}&\|_* \overset{\eqref{ineq:V-cs-1}}{=} \Big\|\sum_{t=0}^k\Big(\Big(W - \frac{1}{N}\mathbf{1}\mathbf{1}^T\Big)\otimes I_m\Big)^{k-t+1}(\rmM^{(t)} - \rmM^{(t-1)})\Big\|_*\nonumber \\
		& \overset{\eqref{def:DeltaM-DeltaV}\eqref{eq:Mk-Mk-1}}{\le} \frac{\theta}{1-\theta}\Big\|\sum_{t=0}^k\Big(\Big(W - \frac{1}{N}\mathbf{1}\mathbf{1}^T\Big)\otimes I_m\Big)^{k-t+1}(\rmG(\rmX^{(t)};\xi^{(t)}) - \nabla \rmF(\rmX^{(t)}))\Big\|_* \nonumber\\ 
		& \qquad + \frac{\theta}{1-\theta}\Big\|\sum_{t=0}^k\Big(\Big(W - \frac{1}{N}\mathbf{1}\mathbf{1}^T\Big)\otimes I_m\Big)^{k-t+1}\Delta \rmM^{(t)}\Big\|_*\nonumber\\
		& \le \frac{\theta}{1-\theta}\Big\|\sum_{t=0}^k\Big(\Big(W - \frac{1}{N}\mathbf{1}\mathbf{1}^T\Big)\otimes I_m\Big)^{k-t+1}(\rmG(\rmX^{(t)};\xi^{(t)}) - \nabla \rmF(\rmX^{(t)}))\Big\|_* + \frac{\theta}{1-\theta}\sum_{t=0}^k\lambda^{k-t+1}\|\Delta \rmM^{(t)}\|_*,\nonumber
	\end{align}
	where the last inequality follows from the triangular inequality and $\|(W - \frac{1}{N}\mathbf{1}\mathbf{1}^T)\otimes I_m\|=\|W - \frac{1}{N}\mathbf{1}\mathbf{1}^T\|\overset{\eqref{upbd:eig-W}}{=}\lambda$. Taking expectation on both sides of this inequality, we obtain that 
	\begin{align}
		\E[\|\Delta \rmV^{(k)}\|_*] & \le \frac{\theta}{1-\theta}\E\Big[\Big\|\sum_{t=0}^k\Big(\Big(W - \frac{1}{N}\mathbf{1}\mathbf{1}^T\Big)\otimes I_m\Big)^{k-t+1}(\rmG(\rmX^{(t)};\xi^{(t)}) - \nabla \rmF(\rmX^{(t)}))\Big\|_*\Big] \nonumber\\ 
		& \qquad + \frac{\theta}{1-\theta}\sum_{t=0}^k\lambda^{k-t+1}\E[\|\Delta \rmM^{(t)}\|_*].   \label{ineq:upbd-Vk-inter-exp}
	\end{align}
	In addition, notice from Assumption \ref{asp:basic}(c) that for each $0\le t\le k$,
	\begin{align}\label{ineq:sdpbd-gradnoise}
		&\E[(\rmG(\rmX^{(t)};\xi^{(t)}) - \nabla \rmF(\rmX^{(t)}))^T(\rmG(\rmX^{(t)};\xi^{(t)}) - \nabla \rmF(\rmX^{(t)}))]\nonumber\\
		&=\sum_{i=1}^N\E[(G_i(X_i^{(t)};\xi_i^{(t)}) - \nabla f_i(X^{(t)}_i))^T(G_i(X_i^{(t)};\xi_i^{(t)}) - \nabla f_i(X^{(t)}_i))] \preceq NV^TV.
	\end{align}
	Hence, Lemma \ref{lem:concentration} holds with $(A_t,U)=(\rmG(\rmX^{(t)};\xi^{(t)}) - \nabla \rmF(\rmX^{(t)}),\sqrt{N}V)$ for $0\le t\le k$. It follows from Lemma \ref{lem:concentration} with $\{(A_t,\Omega_t)\}_{t=0}^k=\{(\rmG(\rmX^{(t)};\xi^{(t)}) - \nabla \rmF(\rmX^{(t)}),((W - \frac{1}{N}\mathbf{1}\mathbf{1}^T)\otimes I_m)^{k-t+1})\}_{t=0}^k$ and $U=\sqrt{N}V$ that
	\begin{align}
		&\E\Big[\Big\|\sum_{t=0}^k\Big(\Big(W - \frac{1}{N}\mathbf{1}\mathbf{1}^T\Big)\otimes I_m\Big)^{k-t+1}(\rmG(\rmX^{(t)};\xi^{(t)}) - \nabla \rmF(\rmX^{(t)}))\Big\|_*\Big]\nonumber\\
		&\le \sqrt{N} \|V\|_*\sqrt{\sum_{t=0}^k\Big\|\Big(W - \frac{1}{N}\mathbf{1}\mathbf{1}^T\Big)\otimes I_m\Big\|^{2(k-t+1)}} = \sqrt{N} \|V\|_*\sqrt{\sum_{t=0}^k\lambda^{2(k-t+1)}}\nonumber\\
		&\le \sqrt{\frac{N\lambda^2}{1-\lambda^2}}\|V\|_*\le \frac{\lambda\sqrt{N}\|V\|_*}{\sqrt{1-\lambda}},\label{ineq:W-inter-concensus}
	\end{align}
	where the last inequality is because $\lambda\in(0,1)$. Substituting this inequality into \eqref{ineq:upbd-Vk-inter-exp}, we obtain that
	\begin{align*}
		\E[\|\Delta \rmV^{(k)}\|_*] \le \frac{\theta\lambda\sqrt{N}\|V\|_*}{(1-\theta)\sqrt{1-\lambda}} + \frac{\theta}{1-\theta}\sum_{t=0}^k\lambda^{k-t+1}\E[\|\Delta \rmM^{(t)}\|_*].    
	\end{align*}
	Summing this inequality over $k=0,\ldots,K-1$, we have
	\begin{align*}
		\sum_{k=0}^{K-1}\E[\|\Delta \rmV^{(k)}\|_*] & \le \frac{K\theta\lambda\sqrt{N}\|V\|_*}{(1-\theta)\sqrt{1-\lambda}} + \frac{\theta}{1-\theta}\sum_{k=0}^{K-1}\sum_{t=0}^k\lambda^{k-t+1}\E[\|\Delta \rmM^{(t)}\|_*]\\
		& \le % \frac{K\theta\lambda\|V\|_*}{1-\theta}\sqrt{\frac{N}{1-\lambda}} + \frac{\theta}{1-\theta}\sum_{t=0}^\infty\lambda^{t+1}\sum_{k=0}^{K-1}\E[\|\Delta M_k\|_*]\\
		\frac{K\theta\lambda\sqrt{N}\|V\|_*}{(1-\theta)\sqrt{1-\lambda}} + \frac{\theta\lambda}{(1-\theta)(1-\lambda)}\sum_{k=0}^{K-1}\E[\|\Delta \rmM^{(k)}\|_*],
	\end{align*}
	where the second inequality is because
	\begin{align*}
		\sum_{k=0}^{K-1}\sum_{t=0}^k\lambda^{k-t+1}\E[\|\Delta \rmM^{(t)}\|_*]\le\bigg(\sum_{k=0}^\infty\lambda^{k+1}\bigg)\cdot\bigg(\sum_{k=0}^{K-1}\E[\|\Delta \rmM^{(k)}\|_*]\bigg) \le \frac{\lambda}{1-\lambda}\sum_{k=0}^{K-1}\E[\|\Delta \rmM^{(k)}\|_*].
	\end{align*}
	Hence, \eqref{ineq:expec-DVk-bd} holds as desired.
\end{proof}

The following lemma provides an estimation error for $\{\rmM^{(k)}\}$ generated by Algorithm \ref{alg:r-msgn-1}.

\begin{lemma}
	Suppose that Assumption~\ref{asp:basic} holds. Let $\{\rmX^{(k)}\}$ be generated by Algorithm~\ref{alg:r-msgn-1} with inputs $(\eta,\theta)$, let $\lambda$, $L_\lambda $, $\{\Delta \rmM^{(k)}\}$, and $V$ be given in \eqref{upbd:eig-W}, \eqref{def:Llambda-objgap}, \eqref{def:DeltaM-DeltaV}, and Assumption \ref{asp:basic}(c), respectively. Then,
	\begin{align}\label{ineq:upbd-exp-DMK}
		\sum_{k=0}^{K-1}\E[\|\Delta \rmM^{(k)}\|_*] \le \frac{1}{\theta}\|\nabla \rmF(\rmX^{(0)})\|_* + \frac{KNL_\lambda \eta}{\theta} + K\sqrt{N\theta}\|V\|_*\qquad K\ge1.
	\end{align}
\end{lemma}

\begin{proof}
	Fix $k\ge0$. It follows from the definition of $\{\Delta \rmM^{(k)}\}$ in \eqref{def:DeltaM-DeltaV} and the update of $\{\rmM^{(k)}\}$ in \eqref{def:ave-update} that 
	\begin{align}
		\Delta \rmM^{(k)} & \overset{\eqref{def:DeltaM-DeltaV}}{=} \rmM^{(k)} - \nabla \rmF(\rmX^{(k)}) \overset{\eqref{def:ave-update}}{=} (1-\theta ) \rmM^{(k-1)} + \theta \rmG(\rmX^{(k)};\xi^{(k)}) - \nabla \rmF(\rmX^{(k)})\nonumber\\
		&=(1-\theta)\Delta \rmM^{(k-1)} + (1-\theta)(\nabla \rmF(\rmX^{(k-1)}) -  \nabla \rmF(\rmX^{(k)}))  + \theta  (\rmG(\rmX^{(k)};\xi^{(k)}) - \nabla \rmF(\rmX^{(k)})).\nonumber
	\end{align}
	Unraveling this recursion for $k+1$ iterations, we obtain that
	\begin{align}
		\Delta \rmM^{(k)} &= (1-\theta)^{k+1} \Delta \rmM^{(-1)} + (1-\theta)\sum_{t=0}^k(1-\theta)^{k-t}(\nabla \rmF(\rmX^{(t-1)}) - \nabla \rmF(\rmX^{(t)})) \nonumber\\
		&\qquad + \theta \sum_{t=0}^k(1-\theta)^{k-t}(\rmG(\rmX^{(t)};\xi^{(t)}) - \nabla \rmF(\rmX^{(t)})).\label{rela:MF-id-1}
	\end{align}
	Recall from \eqref{update-ave-x} that $\|\bar{X}^{(t-1)}-\bar{X}^{(t)}\|=\eta\|\bar{O}^{(t-1)}\|\le\eta$ for all $t\ge1$. By this and the fact $\rmX^{(-1)}=\rmX^{(0)}$, one has that $\|\bar{X}^{(t-1)}-\bar{X}^{(t)}\|\le\eta$ for all $t\ge0$. In addition, by Lemma \ref{lem:tech-do}, Assumption~\ref{asp:basic}(b), and Theorem \ref{thm:cs-error}, one has that for all $t\ge0$,
	\begin{align}
		\|\nabla \rmF(\rmX^{(t-1)}) -  \nabla \rmF(\rmX^{(t)})\|_* & \le \sum_{i=1}^N \|\nabla f_i(X^{(t-1)}_i) - \nabla f_i(X^{(t)}_i) \|_*\nonumber \\
		& \le L_*\sum_{i=1}^N\|X^{(t-1)}_i - X^{(t)}_i\| \le NL_*\|\rmX^{(t-1)} - \rmX^{(t)}\|\nonumber\\
		&\le NL_*(\|\rmX^{(t-1)} - \mathbf{1}\otimes\Bar{X}^{(t-1)}\| + \|\rmX^{(t)} - \mathbf{1}\otimes\Bar{X}^{(t)}\| + \|\Bar{X}^{(t-1)} - \Bar{X}^{(t)}\|)\nonumber \\
		&\le N\bigg(\frac{2\sqrt{N}\lambda}{1-\lambda}+1\bigg)L_*\eta \overset{\eqref{def:Llambda-objgap}}{=} NL_\lambda \eta,\label{ineq:useful-vr}
	\end{align}
	where the first inequality is due to Lemma \ref{lem:tech-do}(ii), the second inequality follows from Assumption \ref{asp:basic}(b), the third inequality follows from Lemma \ref{lem:tech-do}(i), the fourth inequality is due to the triangular inequality, and the fifth inequality follows from Theorem \ref{thm:cs-error} and $\|\Bar{X}^{(t)}-\Bar{X}^{(t-1)}\|\le\eta$. Taking the nuclear norm of \eqref{rela:MF-id-1}, then taking expectations and applying \eqref{ineq:useful-vr}, we obtain that
	\begin{align}
		\E[\|\Delta \rmM^{(k)}\|_*] & \overset{\eqref{rela:MF-id-1}}{\le} (1-\theta)^{k+1}\|\Delta \rmM^{(-1)}\|_* + (1-\theta)\sum_{t=0}^k(1-\theta)^{k-t}\E[\|\nabla \rmF(\rmX^{(t-1)}) - \nabla \rmF(\rmX^{(t)})\|_*] \nonumber\\
		&\qquad + \theta \E\Big[\Big\|\sum_{t=0}^k(1-\theta)^{k-t}(\rmG(\mathbf{X}^{(t)};\xi^{(t)}) - \nabla \rmF(\rmX^{(t)}))\Big\|_*\Big]\nonumber\\
		& \overset{\eqref{ineq:useful-vr}}{\le} (1-\theta)^{k+1}\|\Delta \rmM^{(-1)}\|_* +\frac{NL_\lambda \eta}{\theta} + \theta \E\Big[\Big\|\sum_{t=0}^k(1-\theta)^{k-t}(\rmG(\rmX^{(t)};\xi^{(t)}) - \nabla \rmF(\rmX^{(t)}))\Big\|_*\Big]. \label{ineq:upbd-expec-DeltaM}
	\end{align}
	Recall from \eqref{ineq:sdpbd-gradnoise} that Assumption \ref{asp:basic}(c) implies the following:
	\begin{align*}
		\E[(\rmG(\rmX^{(t)};\xi^{(t)}) - \nabla \rmF(\rmX^{(t)}))^T(\rmG(\rmX^{(t)};\xi^{(t)}) - \nabla \rmF(\rmX^{(t)}))] \preceq NV^TV\qquad\forall t\ge0.
	\end{align*}
	Therefore, Lemma \ref{lem:concentration} holds with $(A_t,U)=(\rmG(\rmX^{(t)};\xi^{(t)}) - \nabla \rmF(\rmX^{(t)}),\sqrt{N}V)$ for any $0\le t\le k$. It then follows from Lemma \ref{lem:concentration} with $\{(A_t,\Omega_t)\}_{t=0}^k=\{(\rmG(\rmX^{(t)};\xi^{(t)}) - \nabla \rmF(\rmX^{(t)}),(1-\theta)^{k-t}I_{Nm})\}_{t=0}^k$ and $U=\sqrt{N}V$ that
	\begin{align}
		\E\Big[\Big\|\sum_{t=0}^k(1-\theta)^{k-t}(\rmG(\rmX^{(t)};\xi^{(t)}) - \nabla \rmF(\rmX^{(t)}))\Big\|_*\Big] & \le \sqrt{N} \|V\|_*\sqrt{\sum_{t=0}^k(1-\theta)^{2(k-t)}} \nonumber \\
        & \le\frac{\sqrt{N}\|V\|_*}{\sqrt{1-(1-\theta)^2}} \le \sqrt{\frac{N}{\theta}}\|V\|_*. \label{ineq:expect-pm-conc}
	\end{align}
	In addition, since we set $\rmM^{(-1)} = \mathbf{0}$ in Algorithm \ref{alg:r-msgn-1} and denote $\rmX^{(-1)} = \rmX^{(0)}$ artificially when defining \eqref{def:DeltaM-DeltaV}, we have $\Delta \rmM^{(-1)} = \nabla \rmF(\rmX^{(0)})$. By this, the above inequality, and \eqref{ineq:upbd-expec-DeltaM}, one has
	\begin{align*}
		\E[\|\Delta \rmM^{(k)}\|_*] \le (1-\theta)^{k+1}\|\nabla \rmF(\rmX^{(0)})\|_*  +\frac{NL_\lambda \eta}{\theta} + \sqrt{N\theta}\|V\|_*.    
	\end{align*}
	Summing this inequality over $k=0,\ldots,K-1$, and using $\theta\in(0,1)$, we obtain that \eqref{ineq:upbd-exp-DMK} holds as desired.
\end{proof}

We are now ready to prove Theorem \ref{thm:stat-error}.

\begin{proof}
	Summing up \eqref{ineq:desc-1} over $k=0,\ldots,K-1$, rearranging the terms, and using \eqref{def:Llambda-objgap} and Assumption \ref{asp:basic}(a), we obtain that
	\begin{align*}
		\frac{1}{K}\sum_{k=0}^{K-1}\|\bar{g}(\rmX^{(k)})\|_* & \le \frac{f(\bar{X}^{(0)}) - f(\bar{X}^{(K)})}{K\eta} + \frac{2}{K}\sum_{k=0}^{K-1}\big(\|\Delta \rmM^{(k)}\|_* + \|\Delta \rmV^{(k)}\|_*\big) + \frac{L_\lambda \eta}{2}  \\
		&\overset{\eqref{def:Llambda-objgap}}{\le} \frac{\Delta_f}{K\eta} + \frac{2}{K}\sum_{k=0}^{K-1}\big(\|\Delta \rmM^{(k)}\|_* + \|\Delta \rmV^{(k)}\|_*\big) + \frac{L_\lambda \eta}{2}.
	\end{align*}
	Taking the expectation on this inequality and using \eqref{ineq:expec-DVk-bd} and \eqref{ineq:upbd-exp-DMK}, we obtain that
	\begin{align*}
		\frac{1}{K}\sum_{k=0}^{K-1}\E[\|\bar{g}(\rmX^{(k)})\|_*] & \le \frac{\Delta_f}{K\eta} + \frac{2}{K}\sum_{k=0}^{K-1}\E[\|\Delta \rmM^{(k)}\|_*] + \frac{2}{K}\sum_{k=0}^{K-1}\E[\|\Delta \rmV^{(k)}\|_*] + \frac{L_\lambda \eta}{2}\\
		&\overset{\eqref{ineq:expec-DVk-bd}}{\le} \frac{\Delta_f}{K\eta} + \frac{2}{K}\bigg[1 + \frac{\theta\lambda}{(1-\theta)(1-\lambda)}\bigg]\sum_{k=0}^{K-1}\E[\|\Delta \rmM^{(k)}\|_*] + \frac{2\theta\lambda\sqrt{N}\|V\|_*}{(1-\theta)\sqrt{1-\lambda}} + \frac{L_\lambda \eta}{2}\\
		&\overset{\eqref{ineq:upbd-exp-DMK}}{\le} \frac{\Delta_f}{K\eta} + \frac{2\theta\lambda\sqrt{N}\|V\|_*}{(1-\theta)\sqrt{1-\lambda}} + \frac{L_\lambda \eta}{2}\\
		&\qquad + \frac{2}{K}\bigg[1 + \frac{\theta\lambda}{(1-\theta)(1-\lambda)}\bigg]\bigg(\frac{1}{\theta}\|\nabla \rmF(\rmX^{(0)})\|_* + \frac{KNL_\lambda \eta}{\theta} + K\sqrt{N\theta}\|V\|_*\bigg) \\
		&\le \frac{\Delta_f}{K\eta} + \frac{2\theta\lambda\sqrt{N}\|V\|_*}{(1-\theta)\sqrt{1-\lambda}} + \frac{(8N+1)L_\lambda \eta}{2\theta(1-\theta)(1-\lambda)} + \frac{4\|\nabla \rmF(\rmX^{(0)})\|_*}{K\theta(1-\theta)(1-\lambda)} + \frac{4\sqrt{N\theta}\|V\|_*}{(1-\theta)(1-\lambda)},
	\end{align*}
	where the last inequality is due to $\theta,\lambda\in(0,1)$. Hence, \eqref{ineq:stat-error} holds as desired.
	
\end{proof}

\subsection{Proof of Theorem \ref{thm:conv}}\label{sec:pf-thm3}

\begin{proof}
	We first prove the first relation of \eqref{upbd:cs-ave}. Notice from the definition of $\hat{\theta}$ in \eqref{def:hat-theta-hat-eta} and $K\ge\frac{4(1-\lambda)\Delta_0L_\lambda}{\|V\|_*^2}$ that $\hat{\theta}\in(0,1/2]$. Using this, and \eqref{ineq:stat-error} with $(\eta,\theta)=(\hat{\eta},\hat{\theta})$, we have
	\begin{align}
		\frac{1}{K} \sum_{k=0}^{K-1} \E[\|\bar{g}(\rmX^{(k)})\|_*] & \overset{\eqref{ineq:stat-error}}{\le} \frac{\Delta_f}{K\hat\eta} + \frac{8\sqrt{N\hat\theta}\|V\|_*}{1-\lambda} + \frac{9NL_\lambda \hat\eta}{\hat\theta(1-\lambda)} + 4\hat\theta\lambda\|V\|_*\sqrt{\frac{N}{1-\lambda}}  + \frac{8\|\nabla \rmF(\rmX^{(0)})\|_*}{K\hat\theta(1-\lambda)}\nonumber \\
		& = 10\sqrt{\frac{N\Delta_fL_\lambda }{(1-\lambda)K\hat{\theta}}} + \frac{8\sqrt{N\hat{\theta}}\|V\|_*}{1-\lambda} + 4\hat\theta\lambda\|V\|_*\sqrt{\frac{N}{1-\lambda}} + \frac{8\|\nabla \rmF(\rmX^{(0)})\|_*}{K\hat\theta(1-\lambda)} \nonumber\\
		&= \frac{18\sqrt{N\|V\|_*}}{(1-\lambda)^{\frac{3}{4}}}\bigg[\frac{\Delta_fL_\lambda }{K}\bigg]^{\frac{1}{4}} + 4\lambda\sqrt{\frac{N\Delta_fL_\lambda }{K}} + \frac{8\|\nabla \rmF(\rmX^{(0)})\|_*\|V\|_*}{(1-\lambda)^{\frac{3}{2}}\sqrt{\Delta_fL_\lambda K}},\nonumber
	\end{align}
	where the first inequality is due to \eqref{ineq:stat-error} and $\hat{\theta}\in(0,1/2]$, the first equality follows from the definition of $\hat{\eta}$ in \eqref{def:hat-theta-hat-eta}, and the last equality follows from the definition of $\hat{\theta}$ in \eqref{def:hat-theta-hat-eta}. Hence, the first relation of \eqref{upbd:cs-ave} holds.
	
	We next prove the second relation of \eqref{upbd:cs-ave}. By \eqref{def:hat-theta-hat-eta}, one has 
	\begin{align*}
		\hat{\eta} = \sqrt{\frac{(1-\lambda)\Delta_f\hat\theta}{NL_\lambda K}} = \bigg[\frac{(1-\lambda)\Delta_f}{K}\bigg]^{\frac{3}{4}}\frac{1}{\sqrt{N\|V\|_*} L_\lambda ^{\frac{1}{4}}}.
	\end{align*}
	Using this and \eqref{ineq:cs-error} with $\eta=\hat{\eta}$, we obtain that
	\begin{align*}
		\|\rmX^{(k)} - \mathbf{1}\otimes\bar{X}^{(k)}\| \le \frac{\sqrt{N}\lambda\hat{\eta}}{1-\lambda} = \Big(\frac{\Delta_f}{K}\Big)^{\frac{3}{4}}\frac{\lambda}{\sqrt{\|V\|_*} [(1-\lambda)L_\lambda ]^{\frac{1}{4}}}.
	\end{align*}
	Hence, \eqref{upbd:cs-ave} holds as desired.
\end{proof}

\section{Proof of the main results in Section \ref{sec:DeMuon+}}\label{apx:pf-DeMuonA}

In this section, we provide the proofs of our main results, particularly, Theorems \ref{thm:stat-error+} and \ref{thm:conv+}.

We define the residual of the $p$th-order Taylor expansion of $\nabla f_i$ as
\begin{align}
R_{p,i}(Y_i,X_i) := \nabla f_i(Y_i) - \sum_{r=1}^p\frac{1}{(r-1)!}\nabla^r f_i(X_i)(Y_i-X_i)^{r-1} \quad \forall X_i,Y_i\in\R^{m\times n},i\in[N].\label{def:Rpi}
\end{align}
We also define its stacked form as
\begin{align}
\rmR_p(\rmY,\rmX) := \left[\begin{matrix}
		R_{p,1}(Y_1,X_1)\\
		\vdots\\
		R_{p,N}(Y_N,X_N)
	\end{matrix}\right]\qquad\forall \rmX,\rmY\in\R^{(Nm)\times n}.\label{def:Rp}
\end{align}

The following lemma provides an upper bound on the residual of the $p$th-order Taylor expansion of $\nabla\rmF$.
\begin{lemma}\label{lem:p-derivative-smooth}
    Suppose that Assumption \ref{asp:high-order} holds. Let $L_{p,*}$ be given in Assumption \ref{asp:high-order}, and let $\|\cdot\|_*$ and $\rmR_p(\cdot,\cdot)$ be defined in \eqref{def:operator-norm} and \eqref{def:Rp}, respectively. Then, it holds that
    \begin{align}\label{ineq:high-order-Lip}
    \|\rmR_p(\rmY,\rmX)\|_*\le\frac{N^pL_{p,*}}{p!}\|\rmY-\rmX\|^p\qquad\forall \rmX,\rmY\in\R^{(Nm)\times n}.    
    \end{align}
\end{lemma}

\begin{proof}
    Fix an arbitrary $U\in\R^{m\times n}$ and $i\in[N]$. Denote $\phi_i(X_i)=\langle\nabla f_i(X_i), U\rangle$. Using this and the definition of $\nabla^{r+1}f_i(X_i)(H)^r$, we obtain that
    \begin{align}\label{relation-Dphi-Dfi}
        \mathcal{D}^r\phi_i(X_i)[V]^r = \langle\nabla^{r+1} f_i(X_i)(V)^r, U\rangle\qquad \forall 1\le r\le p-1, V\in\R^{m\times n}.
    \end{align}
    In addition, using $\phi_i(X_i)=\langle\nabla f_i(X_i), U\rangle$ and \eqref{def:operator-norm}, we have
    \begin{align}\label{ineq:Dphi-Dfi-u}
        \|\mathcal{D}^{p-1}\phi_i(Y_i) - \mathcal{D}^{p-1}\phi_i(X_i)\|_* \le \|U\| \|\mathcal{D}^p f_i(Y_i) - \mathcal{D}^p f_i(X_i)\|_*\qquad\forall X_i,Y_i\in\R^{m\times n}.
    \end{align}
    Fix any $X_i,Y_i\in\R^{m\times n}$, and let $\Delta_i := Y_i-X_i$. By Taylor's expansion, one has that
    \begin{align}
    \phi_i(Y_i) & = \phi_i(X_i) + \sum_{r=1}^{p-2}\frac{1}{r!}\mathcal{D}^r\phi_i(X_i)[\Delta_i]^r + \frac{1}{(p-2)!}\int^1_0(1-t)^{p-2}\mathcal{D}^{p-1}\phi_i(X_i+t\Delta_i)[\Delta_i]^{p-1}\mathrm{d}t\nonumber\\
    & = \phi_i(X_i) + \sum_{r=1}^{p-1}\frac{1}{r!}\mathcal{D}^r\phi_i(X_i)[\Delta_i]^r \nonumber \\
    &\qquad + \frac{1}{(p-2)!}\int^1_0(1-t)^{p-2}(\mathcal{D}^{p-1}\phi_i(X_i+t\Delta_i) - \mathcal{D}^{p-1}\phi_i(X_i))[\Delta_i]^{p-1}\mathrm{d}t. \label{equation:phii-DeltaYX}
    \end{align}
    Using this, \eqref{def:operator-norm}, \eqref{relation-Dphi-Dfi}, and \eqref{ineq:Dphi-Dfi-u}, we obtain that
    \begin{align*}
        & \Big|\Big\langle\nabla f_i(Y_i) - \nabla f_i(X_i) - \sum_{r=1}^{p-1}\nabla^{r+1} f_i(X_i)(\Delta_i)^r, U \Big\rangle\Big| \overset{\eqref{relation-Dphi-Dfi}}{=} \Big|\phi_i(Y_i)-\phi_i(X_i) - \sum_{r=1}^{p-1}\frac{1}{r!}\mathcal{D}^r\phi_i(X_i)[\Delta_i]^r\Big| \\
        & \overset{\eqref{equation:phii-DeltaYX}}{=} \Big|\frac{1}{(p-2)!}\int^1_0(1-t)^{p-2}(\mathcal{D}^{p-1}\phi_i(X_i+t\Delta_i) - \mathcal{D}^{p-1}\phi_i(X_i))[\Delta_i]^{p-1}\mathrm{d}t\Big| \\
        & \overset{\eqref{def:operator-norm}}{\le} \frac{1}{(p-2)!} \|\Delta_i\|^{p-1} \int^1_0 (1-t)^{p-2} \|\mathcal{D}^{p-1}\phi_i(X_i+t\Delta_i) - \mathcal{D}^{p-1}\phi_i(X_i)\|_*\mathrm{d}t\\
        & \overset{\eqref{ineq:Dphi-Dfi-u}}{\le} \frac{1}{(p-2)!} \|\Delta_i\|^{p-1} \|U\| \int^1_0 (1-t)^{p-2} \|\mathcal{D}^{p}f_i(X_i+t\Delta_i) - \mathcal{D}^{p}f_i(X_i)\|_*\mathrm{d}t\\
        & \le \frac{1}{(p-2)!} L_{p,*} \|\Delta_i\|^p \|U\| \int^1_0 (1-t)^{p-2} t\mathrm{d}t = \frac{1}{p!} L_{p,*} \|\Delta_i\|^p \|U\|,
    \end{align*}
    where the last inequality follows from Assumption \ref{asp:high-order}, and the last equality follows from $\int^1_0(1-t)^{p-2}t\mathrm{d}t=1/(p(p-1))$. Taking the maximum of this inequality over all $U$ with $\|U\|\le 1$, and using \eqref{def:Rpi}, we derive that
    \begin{align*}
    \|R_{p,i}(Y_i,X_i)\|_* = \Big\|\nabla f_i(Y_i) - \nabla f_i(X_i) - \sum_{r=1}^{p-1}\nabla^{r+1} f_i(X_i)(\Delta_i)^r\Big\|_* \le \frac{L_{p,*}}{p!} \|Y_i-X_i\|^p.
    \end{align*}
    By this, the definition of $\rmR_p(\cdot,\cdot)$ in \eqref{def:Rp}, and Lemma \ref{lem:tech-do}, one has that
    \begin{align*}
        \|\rmR_{p}(Y_i,X_i)\|_* & \le \sum_{i=1}^N \|R_{p,i}(Y_i,X_i)\|_* \le \frac{L_{p,*}}{p!}\sum_{i=1}^N \|Y_i-X_i\|^p \le \frac{L_{p,*}}{p!} \Big(\sum_{i=1}^N\|Y_i-X_i\|\Big)^p \le \frac{N^pL_{p,*}}{p!} \|\rmY-\rmX\|^p.
    \end{align*}
    Hence, the conclusion of this lemma holds as desired.
\end{proof}

The following lemma is adapted from \citet[Lemma 9]{he2025heavytail}. It provides a set of choices for $\{(\gamma_s,\theta_s)\}$ that satisfy \eqref{linear-cond}.
\begin{lemma}\label{lem:ppt-gamma-theta-choice}
    Let $\gamma\in(0,1/2]$ and a positive integer $q$ be given, and 
    \begin{align*}
        \gamma_s = \frac{\gamma}{s^2},\quad \theta_s = \frac{\prod_{1\le r\le q, r\neq s}(1-r^2/\gamma)}{(s^2/\gamma)\prod_{1\le r\le q,r\neq s}((s^2-r^2)/\gamma)}\qquad\forall 1\le s\le q.
    \end{align*}
    Then, $\{(\gamma_s,\theta_s)\}$ satisfies \eqref{linear-cond}. Moreover, it holds that
    \begin{align}\label{upper-lower-bound-sum-theta}
        \sum_{s=1}^q\theta_s\in\Big(\frac{\gamma}{1+\pi^2/6},2\gamma\Big)\subset(0,1),\quad |\theta_s|\le \frac{4\gamma}{s^2}\qquad \forall s\in[q].
    \end{align}
\end{lemma}

\subsection{Proof of Theorem \ref{thm:stat-error+}}\label{sec:pf-thm2+}

The following lemma provides a descent inequality on the network average $\{\Bar{X}^{(k)}\}$. Its proof is identical to that of Lemma \ref{lem:desc-ave} and is omitted.
\begin{lemma}\label{lem:desc-ave+}
	Suppose that Assumption \ref{asp:basic} holds. Let $\{(\rmX^{(k)},\rmM^{(k)},\rmV^{(k)})\}$ be the sequence generated by Algorithm \ref{alg:r-msgn-1+} with step size $\eta$, and let $L_\lambda $ and $\{(\Delta \rmM^{(k)},\Delta \rmV^{(k)})\}$ be defined in \eqref{def:Llambda-objgap} and \eqref{def:DeltaM-DeltaV}, respectively. Then,
	\begin{align}
		f(\bar{X}^{(k+1)}) \le f(\bar{X}^{(k)}) - \eta\|\bar{g} ({\rmX}^{(k)})\|_* + 2\eta(\|\Delta \rmM^{(k)}\|_* + \|\Delta \rmV^{(k)}\|_*)   + \frac{L_\lambda \eta^2}{2}\qquad\forall k\ge0.\label{ineq:desc-1+}
	\end{align}
\end{lemma}

The following lemma provides a consensus error for $\{\rmV^{(k)}\}$ generated by Algorithm \ref{alg:r-msgn-1+}.
\begin{lemma}
	Suppose that Assumption \ref{asp:basic} and \ref{asp:high-order} hold. Let $\{(\rmM^{(k)},\rmV^{(k)})\}$ be generated by Algorithm \ref{alg:r-msgn-1+} with inputs $\eta$, $q=p-1$, $\{\gamma_s\}$, and $\{\theta_s\}$, let $V$ is given in Assumption \ref{asp:basic}(c), and $\lambda$, $(\theta_{[p]},\theta^\prime_{[p]})$, and $\{(\Delta \rmM^{(k)},\Delta \rmV^{(k)})\}$ be defined as in \eqref{upbd:eig-W}, \eqref{def:Lplambda-sum-theta}, and \eqref{def:DeltaM-DeltaV}, respectively. Then it holds that for all $K\ge1$,
	\begin{align}
		\sum_{k=0}^{K-1}\E[\|\Delta \rmV^{(k)}\|_*] \le \frac{K\theta^\prime_{[p]}\lambda\sqrt{N}\|V\|_*}{(1-\theta_{[p]})\sqrt{1-\lambda}} + \frac{\theta_{[p]}\lambda \sum_{k=0}^{K-1}\E[\|\Delta \rmM^{(k)}\|_*]}{(1-\theta_{[p]})(1-\lambda)}\, + &\, \frac{N\lambda L_\lambda\eta\sum_{s=1}^{p-1}(|\theta_s|/\gamma_s)}{(1-\lambda)(1-\theta_{[p]})}.\label{ineq:expec-DVk-bd+}
	\end{align}
\end{lemma}

\begin{proof}
Fix any $k\ge0$. By the same arguments as for proving \eqref{ineq:V-cs-1}, one has that
\begin{align}
\Delta \rmV^{(k)} = \sum_{t=0}^k\Big(\Big(W - \frac{1}{N}\mathbf{1}\mathbf{1}^T\Big)\otimes I_m\Big)^{k-t+1}\big(\rmM^{(t)} - \rmM^{(t-1)}\big).\label{ineq:V-cs-1+}
\end{align}
Recall from the update of $\{\rmM^{(k)}\}$ in \eqref{def:ave-update+} and the definition of $\theta_{[p]}$ that 
\begin{align*}
\rmM^{(t)} - \rmM^{(t-1)} & = \sum_{s=1}^{p-1} [\theta_s (\rmG(\rmZ^{(t,s)};\xi^{(t,s)}) - \nabla \rmF(\rmZ^{(t,s)}))]  + \theta_{[p]}\cdot (\rmM^{(t)} - \rmM^{(t-1)}) \\
&\qquad + \theta_{[p]}\cdot (\nabla \rmF(\rmX^{(t)}) - \rmM^{(t)}) + \sum_{s=1}^{p-1} [\theta_s (\nabla \rmF(\rmZ^{(t,s)}) - \nabla \rmF(\rmX^{(t)}))] \qquad\forall 0\le t\le k,
\end{align*}
which, by rearranging terms and using $\theta_{[p]}\in(0,1)$, is equivalent to
\begin{align}
\rmM^{(t)} - \rmM^{(t-1)} & = \frac{\sum_{s=1}^{p-1}[\theta_s (\rmG(\rmZ^{(t,s)};\xi^{(t,s)})- \nabla \rmF(\rmZ^{(t,s)}))]}{1-\theta_{[p]}} + \frac{\theta_{[p]}}{1-\theta_{[p]}}\cdot(\nabla \rmF(\rmX^{(t)}) - \rmM^{(t)}) \nonumber\\
&\qquad + \frac{\sum_{s=1}^{p-1}[\theta_s (\nabla \rmF(\rmZ^{(t,s)}) - \nabla \rmF(\rmX^{(t)}))]}{1-\theta_{[p]}} \qquad\forall 0\le t\le k.\label{eq:Mk-Mk-1+}
\end{align}
In addition, by Lemma \ref{lem:tech-do}, Assumption \ref{asp:basic}(b), and Theorem \ref{thm:cs-error+}, one has that for all $t\ge0$,
    \begin{align}
        \|\nabla \rmF(\rmZ^{(t,s)}) - \nabla \rmF(\rmX^{(t)})\|_* & \le \sum_{i=1}^N \|\nabla f_i(Z_i^{(t,s)}) - \nabla f_i(X^{(t)}_i)\|_* \le L_* \sum_{i=1}^N\|Z_i^{(t,s)} - X^{(t)}_i\| \nonumber \\
        &  \le \frac{L_*}{\gamma_s} \sum_{i=1}^N\|X_i^{(t)} - X^{(t-1)}_i\| \le \frac{NL_*}{\gamma_s}\|\rmX^{(t)} - \rmX^{(t-1)}\|  \nonumber \\
        & \le  \frac{N L_*}{\gamma_s} (\|\rmX^{(t-1)} - \mathbf{1}\otimes \Bar{X}^{(t-1)}\| + \|\rmX^{(t)} - \mathbf{1}\otimes \Bar{X}^{(t)}\| + \|\Bar{X}^{(t)} - \Bar{X}^{(t-1)}\|)\nonumber\\
        & \le \frac{N}{\gamma_s}\bigg(\frac{2\sqrt{N}\lambda}{1-\lambda} + 1\bigg) L_*\eta \overset{\eqref{def:Llambda-objgap}}{=} \frac{NL_\lambda \eta}{\gamma_s}.\label{ineq:diff-nabla-ZX}
    \end{align}
	Taking the nuclear norm on both sides of \eqref{ineq:V-cs-1+}, and using the above inequality, \eqref{upbd:eig-W}, \eqref{def:DeltaM-DeltaV}, \eqref{eq:Mk-Mk-1+}, \eqref{ineq:diff-nabla-ZX}, Lemma \ref{lem:tech-do}(i), and Assumption~\ref{asp:basic}(c), we derive that 
	\begin{align}
		\|\Delta \rmV^{(k)}&\|_* \overset{\eqref{ineq:V-cs-1+}}{=} \Big\|\sum_{t=0}^k\Big(\Big(W - \frac{1}{N}\mathbf{1}\mathbf{1}^T\Big)\otimes I_m\Big)^{k-t+1}(\rmM^{(t)} - \rmM^{(t-1)})\Big\|_*\nonumber \\
		& \overset{\eqref{def:DeltaM-DeltaV}\eqref{eq:Mk-Mk-1+}}{\le} \frac{1}{1-\theta_{[p]}}\sum_{s=1}^{p-1}\Big(|\theta_s|\cdot\Big\|\sum_{t=0}^k\Big(\Big(W - \frac{1}{N}\mathbf{1}\mathbf{1}^T\Big)\otimes I_m\Big)^{k-t+1}(\rmG(\rmZ^{(t,s)};\xi^{(t,s)}) - \nabla \rmF(\rmZ^{(t,s)}))\Big\|_*\Big) \nonumber\\ 
		& \qquad + \frac{\theta_{[p]}}{1-\theta_{[p]}}\Big\|\sum_{t=0}^k\Big(\Big(W - \frac{1}{N}\mathbf{1}\mathbf{1}^T\Big)\otimes I_m\Big)^{k-t+1}\Delta \rmM^{(t)}\Big\|_*\nonumber\\
        &\qquad + \frac{1}{1-\theta_{[p]}}\sum_{s=1}^{p-1}\Big(|\theta_s|\cdot\Big\|\sum_{t=0}^k\Big(\Big(W - \frac{1}{N}\mathbf{1}\mathbf{1}^T\Big)\otimes I_m\Big)^{k-t+1}(\nabla \rmF(\rmZ^{(t,s)}) - \nabla \rmF(\rmX^{(t)}))\Big\|_*\Big)\nonumber\\
		& \overset{\eqref{upbd:eig-W}\eqref{ineq:diff-nabla-ZX}}{\le} \frac{1}{1-\theta_{[p]}}\sum_{s=1}^{p-1}\Big(|\theta_s|\cdot\Big\|\sum_{t=0}^k\Big(\Big(W - \frac{1}{N}\mathbf{1}\mathbf{1}^T\Big)\otimes I_m\Big)^{k-t+1}(\rmG(\rmZ^{(t,s)};\xi^{(t,s)}) - \nabla \rmF(\rmZ^{(t,s)}))\Big\|_*\Big)\nonumber \\
        &\qquad + \frac{\theta_{[p]}}{1-\theta_{[p]}}\sum_{t=0}^k(\lambda^{k-t+1}\|\Delta \rmM^{(t)}\|_*) + NL_\lambda\eta\cdot\frac{\sum_{s=1}^{p-1}(|\theta_s|/\gamma_s)}{1-\theta_{[p]}}\cdot\sum_{t=0}^k\lambda^{k-t+1},\nonumber
	\end{align}
	where the last inequality follows from \eqref{ineq:diff-nabla-ZX}, the triangular inequality, and $\|(W - \frac{1}{N}\mathbf{1}\mathbf{1}^T)\otimes I_m\|\overset{\eqref{upbd:eig-W}}{=}\lambda$. Taking expectation on both sides of this inequality, we obtain that 
	\begin{align}
		\E[\|\Delta \rmV^{(k)}\|_*] & \le \frac{1}{1-\theta_{[p]}}\sum_{s=1}^{p-1}\Big(|\theta_s|\cdot\E\Big[\Big\|\sum_{t=0}^k\Big(\Big(W - \frac{1}{N}\mathbf{1}\mathbf{1}^T\Big)\otimes I_m\Big)^{k-t+1}(\rmG(\rmZ^{(t,s)};\xi^{(t,s)}) - \nabla \rmF(\rmZ^{(t,s)}))\Big\|_*\Big]\Big) \nonumber\\ 
		& \qquad + \frac{\theta_{[p]}}{1-\theta_{[p]}}\sum_{t=0}^k(\lambda^{k-t+1}\E[\|\Delta \rmM^{(t)}\|_*]) + NL_\lambda\eta\cdot\frac{\sum_{s=1}^{p-1}(|\theta_s|/\gamma_s)}{1-\theta_{[p]}}\cdot\sum_{t=0}^k\lambda^{k-t+1}.   \label{ineq:upbd-Vk-inter-exp+}
	\end{align}
	By the same arguments as for proving \eqref{ineq:W-inter-concensus}, one has that
	\begin{align*}
		\E\Big[\Big\|\sum_{t=0}^k\Big(\Big(W - \frac{1}{N}\mathbf{1}\mathbf{1}^T\Big)\otimes I_m\Big)^{k-t+1}(\rmG(\rmZ^{(t,s)};\xi^{(t,s)}) - \nabla \rmF(\rmZ^{(t,s)}))\Big\|_*\Big] \le \lambda\|V\|_*\sqrt{\frac{N}{1-\lambda}}.
	\end{align*}
	Substituting this inequality into \eqref{ineq:upbd-Vk-inter-exp+} and using $\sum_{t=1}^k\lambda^t\le\frac{\lambda}{1-\lambda}$ and the definition of $\theta_{[p]}^\prime$ in \eqref{def:Lplambda-sum-theta}, we obtain that
	\begin{align*}
		\E[\|\Delta \rmV^{(k)}\|_*] \le \frac{\theta_{[p]}^\prime\lambda\sqrt{N}\|V\|_*}{(1-\theta_{[p]})\sqrt{1-\lambda}} + \frac{\theta_{[p]}}{1-\theta_{[p]}}\sum_{t=0}^k(\lambda^{k-t+1}\E[\|\Delta \rmM^{(t)}\|_*]) + \frac{N\lambda L_\lambda\eta\sum_{s=1}^{p-1}(|\theta_s|/\gamma_s)}{(1-\lambda)(1-\theta_{[p]})}.
	\end{align*}
	Summing this inequality over $k=0,\ldots,K-1$, we have
	\begin{align*}
		\sum_{k=0}^{K-1}\E[\|\Delta \rmV^{(k)}\|_*] & \le \frac{K\theta_{[p]}^\prime\lambda\sqrt{N}\|V\|_*}{(1-\theta_{[p]})\sqrt{1-\lambda}} + \frac{\theta_{[p]}\sum_{k=0}^{K-1}\sum_{t=0}^k(\lambda^{k-t+1}\E[\|\Delta \rmM^{(t)}\|_*])}{1-\theta_{[p]}}  + \frac{KN\lambda L_\lambda\eta\sum_{s=1}^{p-1}(|\theta_s|/\gamma_s)}{(1-\lambda)(1-\theta_{[p]})}\\
		& \le % \frac{K\theta\lambda\|V\|_*}{1-\theta}\sqrt{\frac{N}{1-\lambda}} + \frac{\theta}{1-\theta}\sum_{t=0}^\infty\lambda^{t+1}\sum_{k=0}^{K-1}\E[\|\Delta M_k\|_*]\\
		\frac{K\theta_{[p]}^\prime\lambda\sqrt{N}\|V\|_*}{(1-\theta_{[p]})\sqrt{1-\lambda}} + \frac{\theta_{[p]}\lambda \sum_{k=0}^{K-1}\E[\|\Delta \rmM^{(k)}\|_*]}{(1-\theta_{[p]})(1-\lambda)}+ \frac{KN\lambda L_\lambda\eta\sum_{s=1}^{p-1}(|\theta_s|/\gamma_s)}{(1-\lambda)(1-\theta_{[p]})},
	\end{align*}
	where the second inequality is because
	\begin{align*}
		\sum_{k=0}^{K-1}\sum_{t=0}^k(\lambda^{k-t+1}\E[\|\Delta \rmM^{(t)}\|_*])\le\bigg(\sum_{k=0}^\infty\lambda^{k+1}\bigg)\cdot\bigg(\sum_{k=0}^{K-1}\E[\|\Delta \rmM^{(k)}\|_*]\bigg) \le \frac{\lambda\sum_{k=0}^{K-1}\E[\|\Delta \rmM^{(k)}\|_*]}{1-\lambda}.
	\end{align*}
	Hence, \eqref{ineq:expec-DVk-bd+} holds as desired.
\end{proof}

The following lemma adapted from \citet[Lemma 10]{he2025heavytail} provides a useful identity.

\begin{lemma}
Suppose that Assumptions \ref{asp:basic} and \ref{asp:high-order} hold. Let $\rmR_p(\cdot,\cdot)$ be defined in \eqref{linear-cond}, and $\{\rmX^{(k)}\}$ and $\{\rmZ^{(k,t)}\}$ be generated by Algorithm \ref{alg:r-msgn-1+} with input parameters $q=p-1$ and $\{(\gamma_s,\theta_s)\}$ satisfying \eqref{linear-cond}, where $p$ is given in Assumption \ref{asp:high-order}. Let $\theta_{[p]}$ be defined in \eqref{def:Lplambda-sum-theta}. Then it holds that for all $k\ge 1$,    
\begin{align}
\nabla\rmF(\rmX^{(k)}) & = (1 - \theta_{[p]}) \nabla\rmF(\rmX^{(k-1)}) + \rmR_p(\rmX^{(k)},\rmX^{(k-1)}) \nonumber \\
&\qquad + \sum_{s=1}^{p-1} (\theta_s\nabla\rmF(\rmZ^{(k,s)})) - \sum_{s=1}^{p-1}(\theta_s\rmR_p(\rmZ^{(k,s)},\rmX^{(k-1)})). \label{eq:identity}
\end{align}

\end{lemma}

The following lemma provides an estimation error for $\{\rmM^{(k)}\}$ generated by Algorithm \ref{alg:r-msgn-1+}.

\begin{lemma}
	Suppose that Assumption~\ref{asp:basic} and \ref{asp:high-order} hold. Let $\{\rmX^{(k)}\}$ be generated by Algorithm~\ref{alg:r-msgn-1+} with input parameters $q=p-1$, $\eta$, $\{(\gamma_s,\theta_s)\}$ satisfying \eqref{linear-cond}, let $\lambda$, $L_\lambda$, $L_{p,\lambda}$, $(\theta_{[p]},\theta_{[p]}^\prime)$, $\{\Delta \rmM^{(k)}\}$, and $V$ be given in \eqref{upbd:eig-W}, \eqref{def:Llambda-objgap}, \eqref{def:Lplambda-sum-theta}, \eqref{def:Lplambda-sum-theta}, \eqref{def:DeltaM-DeltaV}, and Assumption~\ref{asp:basic}(c), respectively. Then it holds that for all $K\ge1$,
	\begin{align}\label{ineq:upbd-exp-DMK+}
		\sum_{k=0}^{K-1}\E[\|\Delta \rmM^{(k)}\|_*] \le \frac{\|\nabla\rmF(\rmX^{(0)})\|_*}{\theta_{[p]}} +\frac{KL_{p,\lambda}\eta^p}{\theta_{[p]}}\bigg(1+\sum_{s=1}^{p-1}\frac{|\theta_s|}{\gamma_s^p}\bigg) +\, &\, K\sqrt{\frac{N}{\theta_{[p]}}}\cdot\theta_{[p]}^\prime\cdot \|V\|_*.
	\end{align}
\end{lemma}

\begin{proof}
	Fix $k\ge0$. It follows from the definition of $\{\Delta \rmM^{(k)}\}$ in \eqref{def:DeltaM-DeltaV} and the update of $\{\rmM^{(k)}\}$ in \eqref{def:ave-update+} that 
	\begin{align}
		\Delta \rmM^{(k)} & \overset{\eqref{def:DeltaM-DeltaV}}{=} \rmM^{(k)} - \nabla \rmF(\rmX^{(k)}) \overset{\eqref{def:ave-update+}}{=} (1-\theta_{[p]}) \rmM^{(k-1)} + \sum_{s=1}^{p-1}(\theta_s\rmG(\rmZ^{(k,s)};\xi^{(k,s)})) - \nabla \rmF(\rmX^{(k)})\nonumber\\
		&\overset{\eqref{eq:identity}}{=}  (1-\theta_{[p]})\Delta \rmM^{(k-1)}  + \sum_{s=1}^{p-1}(\theta_s(\rmG(\rmZ^{(k,s)};\xi^{(k,s)}) - \nabla\rmF(\rmZ^{(k,s)})))\nonumber \\
        &\qquad - \rmR_p(\rmX^{(k)},\rmX^{(k-1)}) + \sum_{s=1}^{p-1}(\theta_s\rmR_p(\rmZ^{(k,s)},\rmX^{(k-1)})).
	\end{align}
	Unraveling this recursion for $k+1$ iterations, we obtain that
	\begin{align}
		\Delta \rmM^{(k)} & = (1-\theta_{[p]})^{k+1}\Delta \rmM^{(-1)} - \sum_{t=0}^k\Big[(1-\theta_{[p]})^{k-t}\Big(\rmR_p(\rmX^{(t)},\rmX^{(t-1)}) - \sum_{s=1}^{p-1}(\theta_s\rmR_p(\rmZ^{(t,s)},\rmX^{(t-1)})) \Big)\Big]\nonumber\\
		&\qquad\quad + \sum_{t=0}^k\Big[(1-\theta_{[p]})^{k-t}\Big(\sum_{s=1}^{p-1}(\theta_s(\rmG(\rmZ^{(t,s)};\xi^{(t,s)}) - \nabla\rmF(\rmZ^{(t,s)})))\Big)\Big].\label{rela:MF-id-1+}
	\end{align}
	Recall from \eqref{update-ave-x} that $\|\bar{X}^{(t-1)}-\bar{X}^{(t)}\|=\eta\|\bar{O}^{(t-1)}\|\le\eta$ for all $t\ge1$. By this and the fact $\rmX^{(-1)}=\rmX^{(0)}$, one has that $\|\bar{X}^{(t-1)}-\bar{X}^{(t)}\|\le\eta$ for all $t\ge0$. It then follows that
    \begin{align}
    \|\rmX^{(t)} - \rmX^{(t-1)}\|^p & \le 3^{p-1}(\|\rmX^{(t-1)} - \mathbf{1}\otimes\Bar{X}^{(t-1)}\|^p + \|\rmX^{(t)} - \mathbf{1}\otimes\Bar{X}^{(t)}\|^p + \|\Bar{X}^{(t-1)} - \Bar{X}^{(t)}\|^p) \nonumber\\
    & \le 3^{p-1}\Big[2\Big(\frac{\sqrt{N}\lambda}{1-\lambda}\Big)^p + 1\Big]\eta^p,\label{ineq:upbd-XtXt-1}
    \end{align}
    where the first inequality is due to $\|A+B+C\|^p\le 3^{p-1}(\|A\|^p + \|B\|^p + \|C\|^p$) for all $p\ge2$ and $A,B,C\in\R^{(Nm)\times n}$, and the second inequality follows from $\|\bar{X}^{(t-1)}-\bar{X}^{(t)}\|\le\eta$ and Theorem \ref{thm:cs-error+}. By this, \eqref{def:Lplambda-sum-theta}, and Lemma \ref{lem:p-derivative-smooth}, one has that for all $t\ge0$,
	\begin{align}
		\|\rmR_p(\rmX^{(t)},\rmX^{(t-1)})\|_* & \le \frac{N^pL_{p,*}}{p!} \|\rmX^{(t)} - \rmX^{(t-1)}\|^p \overset{\eqref{ineq:upbd-XtXt-1}}{\le} \frac{3^{p-1}N^pL_{p,*}}{p!}\Big[2\Big(\frac{\sqrt{N}\lambda}{1-\lambda}\Big)^p + 1\Big]\eta^p\overset{\eqref{def:Lplambda-sum-theta}}{=} L_{p,\lambda} \eta^p,\label{ineq:useful-vr+}\\
        \|\rmR_p(\rmZ^{(t,s)},\rmX^{(t-1)})\|_* & \le \frac{N^pL_{p,*}}{p!} \|\rmZ^{(t,s)} - \rmX^{(t-1)}\|^p \overset{\eqref{def:ave-update+}}{=} \frac{N^pL_{p,*}}{p!\gamma_s^p} \|\rmX^{(t)} - \rmX^{(t-1)}\|^p \nonumber \\
        &\overset{\eqref{ineq:upbd-XtXt-1}}{\le} \frac{3^{p-1}N^pL_{p,*}}{p!}\Big[2\Big(\frac{\sqrt{N}\lambda}{1-\lambda}\Big)^p + 1\Big]\eta^p \overset{\eqref{def:Lplambda-sum-theta}}{=} \frac{L_{p,\lambda}\eta^p}{\gamma_s^p},\label{ineq:useful-vrz+}
	\end{align}
    where the first and third inequalities are due to Lemma \ref{lem:p-derivative-smooth}. Taking the norm of \eqref{rela:MF-id-1+}, then taking expectations and applying \eqref{ineq:useful-vr+} and \eqref{ineq:useful-vrz+}, we obtain that
	\begin{align}
		\E[\|\Delta \rmM^{(k)}\|_*] & \overset{\eqref{rela:MF-id-1+}}{\le} (1-\theta_{[p]})^{k+1}\|\Delta \rmM^{(-1)}\|_* \nonumber\\
        &\qquad + \sum_{t=0}^k\Big[(1-\theta_{[p]})^{k-t}\Big(\E[\|\rmR_p(\rmX^{(t)},\rmX^{(t-1)})\|_*] + \sum_{s=1}^{p-1}(|\theta_s|\cdot\E[\|\rmR_p(\rmZ^{(t,s)},\rmX^{(t-1)}))\|_*] \Big)\Big]\Big] \nonumber\\
		&\qquad + \E\Big[\Big\|\sum_{t=0}^k\Big[(1-\theta_{[p]})^{k-t}\sum_{s=1}^{p-1}(\theta_s(\rmG(\rmZ^{(t,s)};\xi^{(t,s)}) - \nabla\rmF(\rmZ^{(t,s)})))\Big]\Big\|_*\Big] \nonumber\\
		& \overset{\eqref{ineq:useful-vr+}\eqref{ineq:useful-vrz+}}{\le} (1-\theta_{[p]})^{k+1}\|\Delta \rmM^{(-1)}\|_* +\frac{L_{p,\lambda}\eta^p}{\theta_{[p]}}\bigg(1+\sum_{s=1}^{p-1}\frac{|\theta_s|}{\gamma_s^p}\bigg) \nonumber\\
        &\qquad + \sum_{s=1}^{p-1}\Big\{|\theta_s| \cdot \E\Big[\Big\|\sum_{t=0}^k(1-\theta_{[p]})^{k-t}(\rmG(\rmZ^{(t,s)};\xi^{(t,s)}) - \nabla \rmF(\rmZ^{(t,s)}))\Big\|_*\Big]\Big\}. \label{ineq:upbd-expec-DeltaM+}
	\end{align}
	By the same arguments as for proving \eqref{ineq:expect-pm-conc}, one has that
	\begin{align*}
		\E\Big[\Big\|\sum_{t=0}^k(1-\theta_{[p]})^{k-t}(\rmG(\rmZ^{(t,s)};\xi^{(t,s)}) - \nabla \rmF(\rmZ^{(t,s)}))\Big\|_*\Big] \le \sqrt{\frac{N}{\theta_{[p]}}}\|V\|_*.
	\end{align*}
	In addition, since we set $\rmM^{(-1)} = \mathbf{0}$ in Algorithm \ref{alg:r-msgn-1} and denote $\rmX^{(-1)} = \rmX^{(0)}$ artificially when defining \eqref{def:DeltaM-DeltaV}, we have $\Delta \rmM^{(-1)} = \nabla \rmF(\rmX^{(0)})$. By this, the above inequality, and \eqref{ineq:upbd-expec-DeltaM+}, one has
	\begin{align*}
		\E[\|\Delta \rmM^{(k)}\|_*] \le (1-\theta_{[p]})^{k+1}\|\nabla\rmF(\rmX^{(0)})\|_* +\frac{L_{p,\lambda}\eta^p}{\theta_{[p]}}\bigg(1+\sum_{s=1}^{p-1}\frac{|\theta_s|}{\gamma_s^p}\bigg) + \sqrt{\frac{N}{\theta_{[p]}}}\cdot\theta_{[p]}^\prime\cdot\|V\|_*.
	\end{align*}
	Summing this inequality over $k=0,\ldots,K-1$, and using $\theta_{[p]}\in(0,1)$, we obtain that \eqref{ineq:upbd-exp-DMK+} holds.
\end{proof}

We are now ready to prove Theorem \ref{thm:stat-error+}.

\begin{proof}
	Summing up \eqref{ineq:desc-1+} over $k=0,\ldots,K-1$, rearranging the terms, and using \eqref{def:Llambda-objgap} and Assumption \ref{asp:basic}(a), we obtain that
	\begin{align*}
		\frac{1}{K}\sum_{k=0}^{K-1}\|\bar{g}(\rmX^{(k)})\|_* & \le \frac{f(\bar{X}^{(0)}) - f(\bar{X}^{(K)})}{K\eta} + \frac{2}{K}\sum_{k=0}^{K-1}\big(\|\Delta \rmM^{(k)}\|_* + \|\Delta \rmV^{(k)}\|_*\big) + \frac{L_\lambda \eta}{2}  \\
		&\overset{\eqref{def:Llambda-objgap}}{\le} \frac{\Delta_f}{K\eta} + \frac{2}{K}\sum_{k=0}^{K-1}\big(\|\Delta \rmM^{(k)}\|_* + \|\Delta \rmV^{(k)}\|_*\big) + \frac{L_\lambda \eta}{2}.
	\end{align*}
	Taking the expectation on this inequality and using \eqref{ineq:expec-DVk-bd+} and \eqref{ineq:upbd-exp-DMK+}, we obtain that
	\begin{align*}
		& \frac{1}{K}\sum_{k=0}^{K-1}\E[\|\bar{g}(\rmX^{(k)})\|_*] \\
        & \le \frac{\Delta_f}{K\eta} + \frac{2}{K}\sum_{k=0}^{K-1}\E[\|\Delta \rmM^{(k)}\|_*] + \frac{2}{K}\sum_{k=0}^{K-1}\E[\|\Delta \rmV^{(k)}\|_*] + \frac{L_\lambda \eta}{2}\\
		&\overset{\eqref{ineq:expec-DVk-bd+}}{\le} \frac{\Delta_f}{K\eta} + \frac{2\theta_{[p]}^\prime\lambda\sqrt{N}\|V\|_*}{(1-\theta_{[p]})\sqrt{1-\lambda}} + \frac{2N\lambda L_\lambda\eta\sum_{s=1}^{p-1}(|\theta_s|/\gamma_s)}{K(1-\lambda)(1-\theta_{[p]})} + \frac{L_\lambda \eta}{2}  + \frac{4}{K(1-\theta_{[p]})(1-\lambda)} \sum_{k=0}^{K-1}\E[\|\Delta \rmM^{(k)}\|_*]\\
		&\overset{\eqref{ineq:upbd-exp-DMK+}}{\le} \frac{\Delta_f}{K\eta} + \frac{2\theta_{[p]}^\prime\lambda\sqrt{N}\|V\|_*}{(1-\theta_{[p]})\sqrt{1-\lambda}} + \frac{2N\lambda L_\lambda\eta\sum_{s=1}^{p-1}(|\theta_s|/\gamma_s)}{K(1-\lambda)(1-\theta_{[p]})} +\frac{L_\lambda \eta}{2}\\
        &\qquad + \frac{4}{K(1-\theta_{[p]})(1-\lambda)} \bigg[\frac{\|\nabla\rmF(\rmX^{(0)})\|_*}{\theta_{[p]}} +\frac{KL_{p,\lambda}\eta^p}{\theta_{[p]}}\bigg(1+\sum_{s=1}^{p-1}\frac{|\theta_s|}{\gamma_s^p}\bigg) + K\sqrt{\frac{N}{\theta_{[p]}}}\cdot\theta_{[p]}^\prime\cdot\|V\|_*\bigg]\\
		&\le \frac{\Delta_f}{K\eta} + \frac{2\theta_{[p]}^\prime\lambda\sqrt{N}\|V\|_*}{(1-\theta_{[p]})\sqrt{1-\lambda}} + \frac{2N\lambda L_\lambda\eta\sum_{s=1}^{p-1}(|\theta_s|/\gamma_s)}{K(1-\lambda)(1-\theta_{[p]})} + \frac{L_\lambda\eta}{2} + \frac{4\|\nabla\rmF(\rmX^{(0)})\|_*}{K\theta_{[p]}(1-\theta_{[p]})(1-\lambda)}\\
        &\qquad +\frac{4L_{p,\lambda}\eta^p(1+\sum_{s=1}^{p-1}(|\theta_s|/\gamma_s^p))}{\theta_{[p]}(1-\theta_{[p]})(1-\lambda)} + \frac{4\sqrt{N}\theta_{[p]}^\prime\|V\|_*}{\sqrt{\theta_{[p]}}(1-\theta_{[p]})(1-\lambda)},
	\end{align*}
	where the last inequality is due to $\theta_{[p]}\in(0,1)$ and $\lambda\in(0,1)$. Hence, \eqref{ineq:stat-error+} holds as desired.
	
\end{proof}

\subsection{Proof of Theorem \ref{thm:conv+}}\label{subsec:proof-thm-conv+}

\begin{proof}
We first prove the first relation in \eqref{upbd:cs-ave+}. Let $\hat{\theta}_{[p]}=\sum_{r=1}^{p-1}\hat{\theta}_s$ and $\hat{\theta}_{[p]}^\prime=\sum_{r=1}^{p-1}|\hat{\theta}_s|$. Recall from \eqref{def:hat-theta-hat-eta++} and Lemma \ref{lem:ppt-gamma-theta-choice} that $\hat{\theta}_{[p]}\in(\hat{\gamma}/3,2\hat{\gamma})$ and $|\hat{\theta}_s|\le4\hat{\gamma}/s^2$ for all $1\le s\le p-1$. By $\hat{\theta}_{[p]}<2\hat{\gamma}$, the definition of $\hat{\gamma}$, and $K\ge\frac{16L_{p,\lambda}^{1/p}(1-\lambda)\Delta_f}{(\sqrt{N}\|V\|_*)^{(p+1)/p}}$, one has $\hat{\theta}_{[p]}\in(0,1/2)$. Also, we have $\hat{\theta}_{[p]}^\prime\le4\hat{\gamma}\sum_{s=1}^{p-1}1/s^2\le4\hat{\gamma}\sum_{s=1}^\infty1/s^2 = 2\pi^2\hat{\gamma}/3<7\hat{\gamma}$. Using these and \eqref{ineq:stat-error+} with $(\eta,\theta_{[p]},\theta_{[p]}^\prime)=(\hat{\eta},\hat\theta_{[p]},\hat\theta_{[p]}^\prime)$, and $\{(\gamma_s,\theta_s)\}=\{(\hat\gamma_s,\hat{\theta}_s)\}$, we obtain that 
\begin{align*}
\frac{1}{K} \sum_{k=0}^{K-1} \E[\|\bar{g}(\rmX^{(k)})\|_*] &\le \frac{\Delta_f}{K\hat\eta} +\frac{8L_{p,\lambda}\hat\eta^p(1+\sum_{s=1}^{p-1}(|\hat\theta_s|/\hat\gamma_s^p))}{\hat\theta_{[p]}(1-\lambda)} + \frac{8\sqrt{N}\hat\theta_{[p]}^\prime\|V\|_*}{\sqrt{\hat\theta_{[p]}}(1-\lambda)} + \frac{L_\lambda\hat\eta}{2} \nonumber\\
&\qquad + \frac{4\hat\theta_{[p]}^\prime\lambda\sqrt{N}\|V\|_*}{\sqrt{1-\lambda}} + \frac{4N\lambda L_\lambda\hat\eta\sum_{s=1}^{p-1}(|\hat\theta_s|/\hat\gamma_s)}{K(1-\lambda)} + \frac{8\|\nabla\rmF(\rmX^{(0)})\|_*}{K\hat\theta_{[p]}(1-\lambda)}\\
& \le \frac{\Delta_f}{K\hat\eta} +\frac{24L_{p,\lambda}\hat\eta^p(1+(4/\hat{\gamma}^{p-1})\sum_{s=1}^{p-1}s^{2(p-1)})}{\hat{\gamma}(1-\lambda)} + \frac{97\sqrt{N\hat{\gamma}}\|V\|_*}{1-\lambda} + \frac{L_\lambda\hat\eta}{2} \\
& \qquad + \frac{28\hat\gamma\lambda\sqrt{N}\|V\|_*}{\sqrt{1-\lambda}} + \frac{16(p-1)N\lambda L_\lambda\hat\eta}{K(1-\lambda)} + \frac{24\|\nabla\rmF(\rmX^{(0)})\|_*}{K\hat{\gamma}(1-\lambda)} \\
& \le \frac{\Delta_f}{K\hat\eta} +\frac{96(p-1)^{2p-1}L_{p,\lambda}\hat\eta^p}{(1-\lambda)\hat{\gamma}^p} +  \frac{24L_{p,\lambda}\hat\eta^p}{(1-\lambda)\hat{\gamma}} + \frac{97\sqrt{N\hat{\gamma}}\|V\|_*}{1-\lambda} + \frac{L_\lambda\hat\eta}{2}\\
& \qquad + \frac{28\hat\gamma\lambda\sqrt{N}\|V\|_*}{\sqrt{1-\lambda}} + \frac{16(p-1)N\lambda L_\lambda\hat\eta}{K(1-\lambda)} + \frac{24\|\nabla\rmF(\rmX^{(0)})\|_*}{K\hat{\gamma}(1-\lambda)} \\
& = \mathcal{O}(K^{-p/(3p+1)}).
\end{align*}
Hence, the first relation of \eqref{upbd:cs-ave+} holds.

We next prove the second relation of \eqref{upbd:cs-ave+}. It follows from Theorem \ref{thm:cs-error+} and \eqref{def:hat-theta-hat-eta+} that
\begin{align*}
    \max_{0\le k\le K-1}\{\|\rmX^{(k)} - \mathbf{1}\otimes\bar{X}^{(k)}\|\} \le \mathcal{O}(\hat\eta) = \mathcal{O}(K^{-(2p+1)/(3p+1)}).
\end{align*}
Hence, the second relation of \eqref{upbd:cs-ave+} holds as desired.
    
\end{proof}

\section{Extra experimental setups}
\label{app:params}

\begin{table*}[htpb]
\centering
\caption{Parallel update schemes and hyperparameters for competing methods, {\it DSGD} \citep{nedic2009distributed}, {\it DSGD-C} \citep{dsgd-clip}, {\it DSGD-N} \citep{yu2026decentralized} and our {\it DeMuon}. For each method and topology, the learning rate $\eta$ and decay schedule are chosen to minimize consensus error, enabling a fair comparison of training quality.}
\label{tab:stepsizes}
\small
\begin{tabular}{p{1.7cm}| p{6.5cm} p{4.9cm}}
\hline
{\it Algorithm} & Parallel Update Scheme at Agent $i$  & {Hyperparameters} \\
\hline
&& \\
\multirow{3}{*}{\it DSGD} & \multirow{3}{*}{\makecell[l]{$X_i^{k+1}=\sum_{j=1}^Nw_{ij}(X_j^k-\eta_k G(X_j^k;\xi_j^k))$}} & Complete: $\eta_k=0.006/\sqrt{k}$\\
&& Exp: $\eta_k=0.03(1-k/T)$\\
&& Ring: $\eta_k=0.03(1-k/T)$\\
&& \\
\hline
&& \\
\multirow{3}{*}{\it DSGD-C} & \multirow{3}{*}{\makecell[l]{$\tau_k=\tau k^{2/5}$\\$X_i^{k+1}=\sum_{j=1}^Nw_{ij}(X_j^k-\eta_k \text{clip}(G(X_j^k;\xi_j^k),\tau_k))$}} & Complete: $\eta_k=0.6/k,\ \tau=0.1$\\
&& Directed Exp.: $\eta_k=0.2(1-k/T),\ \tau=0.1$\\
&& Ring: $\eta_k=0.1(1-k/T),\ \tau=0.1$\\
&& \\
\hline
&& \\
\multirow{3}{*}{\it DSGD-N} & \multirow{3}{*}{\makecell[l]{$M_i^k = (1-\theta )M_i^{k-1} +\theta  G(X^k_i;\xi^k_i)$  \\
$V_i^k = \sum_{j=1}^N w_{ij} (V_j^{k-1} + M_j^k - M_j^{k-1})$\\
$X_i^{k+1}=\sum_{j=1}^N w_{ij}(X_j^k-\eta_k V_j^k/\|V_j^k\|_F)$}} & Complete: $\eta_k=0.07,\ \theta=0.2$\\
&& Directed Exp.: $\eta_k=0.05(1-k/T),\ \theta=0.2$\\
&& Ring: $\eta_k=0.03(1-k/T),\ \theta=0.2$\\
&& \\
\hline
&& \\
\multirow{3}{*}{{\it DeMuon}} & \multirow{3}{*}{\makecell[l]{See steps \eqref{update-mk}-\eqref{update-xk}}} & Complete: $\eta_k=0.1/\sqrt{k},\ \theta=0.8$\\
&& Directed Exp.: $\eta_k=0.005(1-k/T),\ \theta=0.2$\\
&& Ring: $\eta_k=0.003(1-k/T),\ \theta=0.2$\\
&& \\
\hline
\end{tabular}
\end{table*}

\end{document}